\documentclass[11pt]{article}
\usepackage[margin=1in]{geometry}

\usepackage[pdftex]{graphicx}
\graphicspath{{../pdf/}{../jpeg/}}
\DeclareGraphicsExtensions{.pdf,.jpeg,.png}
\usepackage{amssymb,amsmath,url}
\usepackage{bm}

\usepackage{multirow}
\usepackage{booktabs}	
\usepackage{epstopdf}
\usepackage{color}
\usepackage{subfigure}
\usepackage{bigfoot}		
\usepackage{algpseudocode}
\usepackage{pifont}
\usepackage{tikz}
\usepackage{color}
\usepackage{cite}
\usepackage{lipsum}
\usepackage{algorithm} 
\usepackage[pdftex]{graphicx}

\algrenewcommand\algorithmicrequire{\textbf{Input:}}
\algrenewcommand\algorithmicensure{\textbf{Output:}}
\algrenewcommand\algorithmicforall{\textbf{For}}

\newtheorem{theorem}{Theorem}
\newtheorem{lemma}{Lemma}
\newtheorem{proposition}{Proposition}

\allowdisplaybreaks
\newtheorem{definition}{Definition}

\newtheorem{remark}{Remark}

\newtheorem{assumption}{Assumption}
\DeclareMathOperator*{\argmin}{arg\,min}
\DeclareMathOperator*{\argmax}{arg\,max}

\title{\LARGE
	Safe Zeroth-Order Optimization Using Quadratic Local Approximations
}
\author{Baiwei Guo\thanks{B. Guo and G. Ferrari-Trecate are with the DECODE group, Institute of Mechanical Engineering, EPFL, Switzerland. e-mails: {\tt \{baiwei.guo, giancarlo.ferraritrecate\}@epfl.ch}.} \and Yuning Jiang \thanks{Y. Jiang is with the PREDICT group, Institute of Mechanical Engineering, EPFL, Switzerland. e-mails: {\tt yuning.jiang@epfl.ch}.} \and Giancarlo Ferrari-Trecate \footnotemark[1] \and Maryam Kamgarpour\thanks{ M. Kamgarpour is with the SYCAMORE group, Institute of Mechanical Engineering, EPFL, Switzerland. e-mails: {\tt maryam.kamgarpour@epfl.ch}.}}

\allowdisplaybreaks
\begin{document}
\maketitle \let\thefootnote\relax\footnotetext{This work was supported by the Swiss National Science Foundation under the NCCR Automation (grant agreement 51NF40\textunderscore 80545).}\\

\begin{abstract}\noindent {This paper addresses smooth constrained optimization problems with unknown objective and constraint functions.} The main goal of this work is to generate a sequence of feasible (herein, referred to as safe) points converging towards a KKT primal-dual pair. Assuming to have prior knowledge on the smoothness of the unknown functions, we propose a novel zeroth-order method that iteratively computes quadratic approximations of the constraint functions, constructs local feasible sets, and optimizes over them. We prove that this method returns an $\eta$-KKT pair within $O({d}/{\eta^{2}})$ iterations {and $O({d^2}/{\eta^{2}})$ samples (where $d$ is the problem dimension) while every sample is within the feasible set}. Moreover, we numerically show that our method can achieve fast convergence compared with some state-of-the-art zeroth-order safe approaches. The effectiveness of the proposed approach is also illustrated by applying it to nonconvex optimization problems in optimal control and power system operation.
\end{abstract}
\section{Introduction}
\label{sec: introduction}
Applications ranging from power network operations \cite{chu2021frequency}, machine learning \cite{chen2020DFO} and trajectory optimization \cite{manchester2016derivative} to optimal control~\cite{rao1998application,xu2022vabo} require solving complex optimization problems where feasibility (i.e., the fulfillment of the hard constraints) is essential. However, in practice, we do not always have access to the expressions of the objective and constraint functions. 

To address an unmodeled constrained optimization, we develop a safe zeroth-order optimization method in this paper. Zeroth-order methods rely only on sampling (i.e., evaluating the unknown objective and constraint functions at a set of chosen points) \cite{Bajaj2021}. Safety, herein referring to the feasibility of the samples, is essential in several real-world problems, e.g., medical applications \cite{yousefi2017formally}  and racing car control \cite{hewing2019cautious}. Below, we review the pertinent literature on zeroth-order optimization, highlighting, specifically, safe zeroth-order methods. 

Classical techniques for zeroth-order optimization can be classified as direct-search-based (where a set of points around the current point is searched for a lower value of the objective function), gradient-descent-based (where the gradients are estimated based on samples), and model-based (where a local model of the objective function around the current point is built and used for local optimization)~\cite[Chapter~9]{jorge2006numerical}. Examples of these three categories for unconstrained optimization are, respectively, pattern search methods \cite{lewis2000pattern}, randomized stochastic gradient-free methods \cite{ghadimi2013stochastic}, and trust region methods \cite{conn2000trust}. Pattern search methods are extended in \cite{lewis2002globally} to solve optimization problems with known constraints. 

In case the explicit formulations of both objective and constraint functions are not available, {the work \cite{amani2019linear} solves the problem by learning the functions using non-parametric models. However, this method only addresses linear programs.} When the unmodelled constraints are nonlinear, one can use two-phase methods \cite{echebest2017inexact,bajaj2018trust} where an optimization phase reduces the objective function subject to relaxed constraints and a restoration phase modifies the result of the first phase to regain feasibility. A drawback of these approaches is the lack of a guarantee for sample feasibility (i.e., each sample satisfying the constraints). {Therefore, they cannot be used for optimization tasks with hardware in the loop, since any infeasible sample may damage the hardware.}

For sample feasibility, the zeroth-order methods of \cite{sui2015safe,turchetta2019safe,sabug2022smgo}, including SafeOpt and its variants, assume the knowledge of a Lipschitz constant of the objective and constraint functions, while \cite{vinod2022constrained} utilizes the Lipschitz constants of the gradients of these functions (the smoothness constants). With these quantities, one can build local proxies for the constraint functions and provide a confidence interval for the true function values. By starting from a feasible point, \cite{sui2015safe,sabug2022smgo,vinod2022constrained} utilize the proxies to search for potential minimizers. However, for each search, one may have to use a global optimization method to solve a non-convex problem, which makes the algorithm computationally intractable if there are many decision variables.

Another research direction aimed at the feasibility of the samples is to include barrier functions in the objective to penalize the proximity to the boundary of the feasible set \cite{lewis2002globally,audet2009progressive}. In this category, extremum seeking methods estimate the gradient of the new objective function by adding sinusoidal perturbations to the decision variables~\cite{arnold2015model}. However, due to the perturbations, these methods have to adopt a sufficiently large penalty coefficient to ensure all the samples fall in the feasible region. This strategy sacrifices optimality since deriving a near-optimal solution requires a small penalty coefficient.  In contrast, the LB-SGD algorithm proposed in~\cite{usmanova2022log} uses log-barrier functions and ensures the feasibility of the samples despite a small penalty coefficient. After calculating a descent direction for the cost function with log-barrier penalties, this method exploits the Lipshcitz and smoothness constants of the constraint functions to build local safe sets for selecting the step size of the descent. Although LB-SGD comes with a polynomial worst-case complexity in problem dimension, it might converge slowly, {\color{black}even for convex problems}. The reason is that as the iterates approach the boundary of the feasible set, the log-barrier function and its derivative become very large, leading to very conservative local feasible sets and slow progress of the iterates. 

Safe zeroth-order optimization has been an increasingly important topic in the learning-based control community. One application is constrained optimal controller tuning with unknown system dynamics. In reinforcement learning, Constrained Policy Optimization \cite{achiam2017constrained} and Learning Control Barrier Functions \cite{qin2022sablas} (model-free) are used to find the optimal safe controller, but feasibility during training cannot be ensured. Bayesian Optimization can also be applied to optimal control in a zeroth-order manner. For example, \cite{xu2022vabo} proposes Violation-Aware Bayesian Optimization to optimize three set points of a vapor compression system, \cite{berkenkamp2016safe} utilizes SafeOpt to tune a linear control law with two parameters for quadrotors, and \cite{konig2021safe} implements the Goal Oriented Safe Exploration algorithm in \cite{turchetta2019safe} to optimize a PID controller with three parameters for a rotational axis drive. Although these variants of Bayesian Optimization offer guarantees of sample feasibility, they scale poorly to high-dimensional systems due to the non-convexity of the subproblems and the need for numerous samples.

\textit{Contributions:} We  develop a zeroth-order method for smooth optimization problems with guaranteed sample feasibility and convergence. The approach is based on designing quadratic local proxies of the objective and constraint functions. Preliminary results, presented in \cite{guo2022safe}, focused on the case of convex objective and constraints. There, we proposed an algorithm that sequentially solves convex Quadratically Constrained Quadratic Programming (QCQP) subproblems. We showed that all the samples are feasible, and one accumulation point of the iterates is the minimizer. This paper significantly extends \cite{guo2022safe} in the following aspects:
\begin{enumerate}
\item[(a)] we show in Section \ref{sec: convergence of x_K} that, under mild assumptions, our safe zeroth-order algorithm has iterates whose accumulation points are KKT pairs even for \textbf{non-convex} problems;
\item[(b)] besides the asymptotic results in (a), given $\eta>0$, we add \textbf{termination conditions} to the zeroth-order algorithm and guarantee in Section \ref{sec: eta-KKT} that the returned primal-dual pair is an \textbf{$\eta$-KKT} pair (see Definition \ref{def: KKT}) of the optimization problem. {We further show in Section \ref{sec: complexity_analysis} that under mild assumptions our algorithm terminates in $O(\frac{d}{\eta^2})$ iterations and requires $O(\frac{d^2}{\eta^2})$ samples where $d$ is the number of the decision variables.}
\item[(c)] we present in Section \ref{sec: experiments} an example illustrating that our algorithm achieves faster convergence than state-of-the-art zeroth-order methods that guarantee sample feasibility. We further apply the algorithm to optimal control and optimal power flow problems, showing that the results returned by our algorithm are almost identical to those provided by commercial solvers utilizing the true model.
\end{enumerate}
\textit{Notations:}
We use $e_i\in\mathbb R^d$ to define the $i$-th standard basis of the vector space $\mathbb{R}^d$ and $\|\cdot\|$ to denote the 2-norm throughout the paper. Given a vector $x\in \mathbb{R}^d$ and a scalar $\epsilon>0$, we write $x = [x^{(1)}, \ldots, x^{(d)}]^\top$, $\mathcal B_\epsilon(x)=\{y:\|y-x\|\leq \epsilon\}$ and $\mathcal{SP}_\epsilon(x)=\{y:\|y-x\|= \epsilon\}$. {We use $\mathbb{Z}_i^j=\{i,i+1,\ldots,j\}$ to define the set of integers ranging from $i$ to $j$ with $i<j$. We also define $\mathbf{I}_n$ as the identity matrix with $n$ rows. Given a list of vectors $a_1,a_2,\ldots,a_n\in\mathbb{R}^d$, we use $\{a_k\}^n_{k=1}$ to denote the set $\{a_1,a_2,\ldots,a_n\}$.}

\section{Problem Formulation}
\label{sec::Prob}
We address the constrained optimization problem in the form
\begin{equation}
\label{eq: optimization problem}
\min_{x\in \mathbb{R}^d}  \;\; f_0(x)\quad \text{subject to}\;\;
x\in \Omega
\end{equation}
with feasible set $\Omega:=\{x\in \mathbb{R}^d: f_i(x)\leq 0, i\in\mathbb Z_1^m\}$.
We consider the setting where the continuously differentiable functions $f_i:\mathbb R^d\to\mathbb R$ are not explicitly known but can be queried. {If there is hardware in the loop, each query (also called sample) requires an experiment and measurements on the physical system. Therefore, sample feasibility is essential. In this paper, we aim to derive a local optimization algorithm that generates only feasible samples and for any given $\eta>0$ returns an $\eta$-approximate KKT pair of \eqref{eq: optimization problem} defined as follows.}
\begin{definition}
\label{def: KKT}
For $\eta>0$, a pair $(x,\lambda)$ with $x\in \Omega$ and $\lambda\in \mathbb{R}^m_{\geq 0}$ is an $\eta$-approximate KKT ($\eta$-KKT for short) pair of the problem \eqref{eq: optimization problem} if 
\begin{subequations}
\label{eq:kkt}
\begin{align}
\|\nabla f_0(x) + \sum^m_{i=1}\lambda^{(i)}\nabla f_i(x)\| &\leq \eta,\\
|\lambda^{(i)}f_i(x)| &\leq \eta,\quad i\in \mathbb{Z}^{m}_{1}.
\end{align}
\end{subequations}
If $(x^*,\lambda^*)\in \Omega\times \mathbb{R}^m_{\geq 0}$ fulfills~\eqref{eq:kkt} with $\eta = 0$, we say that it is a KKT pair.
\end{definition}
For any optimization problem with differentiable objective and
constraint functions for which strong duality holds, any pair of primal and dual optimal points must be a KKT pair. If the optimization problem is furthermore convex, any KKT pair satisfies primal and dual optimality \cite{boyd2004convex}. The optimization methods that aim to obtain a KKT pair, such as Newton-Raphson and interior point methods, might converge to a local optimum. Despite this drawback, these local methods are extensively applied, because local algorithms are more efficient to implement and KKT pairs are good enough for some applications, such as machine learning \cite{jain2017non}, optimal control \cite{rawlings2017model} and optimal power flow \cite{kundur2022power}. Considering that, in general, numerical solvers cannot return an exact KKT pair, the concept of $\eta$-KKT pair indicates how close primal and dual solutions are to a KKT pair \cite{dutta2013approximate}. According to \cite[Theorem 3.6]{dutta2013approximate}, under mild assumptions, one can make the $\eta$-KKT pair arbitrarily close (in Euclidean distance) to an KKT pair of \eqref{eq: optimization problem} by decreasing $\eta$. In many numerical optimization methods \cite{grapiglia2022tensor,dussault2020approximate}, one can trade off accuracy against efficiency by tuning $\eta$.

We assume, without loss of generality, the objective function $f_0(x)$ is explicitly known and linear. Indeed, when the function $f_0(x)$ in \eqref{eq: optimization problem} is not known but can be queried, the problem in \eqref{eq: optimization problem} can be written as 
\begin{align*}
\min_{(x,\gamma)\in \mathbb{R}^{d+1}
} \quad & \gamma\\
\mathrm{subject\;to}  \quad & f_0(x)-\gamma\leq 0,\\
& f_i(x)\leq 0, \quad i \in\mathbb Z_1^m,
\end{align*}
where the objective function is now known and linear. 
Throughout this paper, we make the following assumptions on the smoothness of the objective and constraint functions and the availability of a strictly feasible point.
\begin{assumption}
\label{ass: smoothness}
The functions $f_i(x)$, $i\in\mathbb Z_0^m$ are continuously differentiable and {we know constants $L_i,M_i>0$ and a set $\overline{\Omega}$ such that $\Omega\subset \overline{\Omega}$ and for any $x_1$, $x_2\in \overline{\Omega}$,}
\begin{subequations}
\label{eq: smoothness}
\begin{align}
\label{eq: smoothness_1}
|f_i(x_1)-f_i(x_2)|&\leq L_i\|x_1-x_2\|,\\
\label{eq: smoothness_2}
\|\nabla f_i(x_1)-\nabla f_i(x_2)\|&\leq M_i\|x_1-x_2\|.
\end{align}
\end{subequations}
We also assume that the known Lipschitz and smoothness constants ${L}_i$ and ${M}_i$ verify that 
\begin{equation}
\label{eq: LANDM}
{L}_i>L_{i,\mathrm{inf}} \text{ and }{M}_i> M_{i,\mathrm{inf}},
\end{equation}
where \begin{align*}
L_{i,\mathrm{inf}}:=&\inf\{L_i: \text{ }\eqref{eq: smoothness_1} \text{ holds},\forall x_1,x_2\in \Omega\},\\
M_{i,\mathrm{inf}}:=&\inf\{M_i: \text{ }\eqref{eq: smoothness_2} \text{ holds},\forall x_1,x_2\in \Omega\}.
\end{align*}
\end{assumption}
In the remainders of this paper, we also define $L_{\max} = \max_{i\geq 1}L_i$ and  $M_{\max} = \max_{i\geq 1}M_i$.
\begin{remark}
{While the set $\Omega$ is unknown, one can define the set $\overline{\Omega}$ by relying on prior knowledge. For example, in optimal power flow problems (see Section 6.3), the decision variables are generator voltage and active power output. The ranges of these variables are part of the generator specifications.} 
\vspace{-0.2cm}

The bounds in \eqref{eq: smoothness} are utilized in several works on zeroth-order optimization, e.g., \cite{wibisono2012finite,conn2009introduction}. As it will be clear in the sequel, these bounds allow one to estimate the error of local approximations of the unknown functions and their derivatives. In practice, it is usually impossible to obtain $L_{i,\mathrm{inf}}$ and $M
_{i,\mathrm{inf}}$, thus we only assume to know the upperbounds ${L}_i>L_{i,\mathrm{inf}}$ and ${M}_i> M_{i,\mathrm{inf}}$. In case ${L}_i$ and ${M}_i$ are not known, we regard them as hyperparameters and describe how to tune them in Remark \ref{rmk: LM}.
\end{remark}
\begin{assumption}
\label{ass: strict feasible}
There exists a known strictly feasible point $x_0$, i.e., $f_i(x_0)<0$ for all $i\in\mathbb Z_1^m$. 
\end{assumption}
\begin{remark}
{For convex optimization problems, the existence of a strictly feasible point is called Slater’s Condition and is commonly assumed in several optimization methods \cite[Section 5.3.2]{boyd2004convex}.} Moreover, several works on safe learning \cite{sui2015safe,usmanova2022log} assume a strictly feasible point used for initializing the algorithm. Assumption \ref{ass: strict feasible} is necessary for designing an algorithm whose iterates remain feasible since the constraint functions are unknown a priori. Practically, it holds in several applications. For example, in any robot mission planning,  the robot is placed initially at a safe point and needs to gradually explore the neighboring regions while ensuring the feasibility of its trajectory. Similarly, in the optimization of manufacturing processes, often an initial set of (suboptimal) design parameters satisfying the problem constraints are known \cite{rupenyan2021performance}. Another example is frequency control of power grids, where the initial frequency is guaranteed to lie within certain bounds by suitably regulating the power reserves and loads \cite{kundur2022power}. 
\end{remark}
\begin{assumption}
\label{ass: bounded_sublevel}
There exists $\beta\in \mathbb{R}$ such that the sublevel set $\mathcal{P}_\beta=\{x\in \Omega:f_0(x)\leq \beta\}$ is bounded and includes the initial feasible point $x_0$.
\end{assumption}
Under Assumption \ref{ass: bounded_sublevel}, for any iterative algorithm producing non-increasing objective function values $\{f_0(x_k)\}_{k\geq 0}$, we ensure the iterates $\{x_{k}\}_{k\geq 0}$ to be within the bounded set $\mathcal{P}_\beta$.

We highlight that Assumptions 1-3 stand \textit{throughout this paper}. By exploiting them, we design in the following section an algorithm that iteratively optimizes $f_0(x)$. 

\section{Safe Zeroth-Order Algorithm}
\label{sec: algorithm_tot}
Before introducing the iterative algorithm, this section proposes an approach to construct local feasible sets by using samples around a given strictly feasible point. To do so, we recall the properties of a gradient estimator constructed through finite differences.

The gradients of the unknown functions $\{f_i\}_{i=1}^m$ can be approximated as  
\begin{equation}
\label{eq: gradient approximation}
{\nabla}^\nu f_i\left(x\right):=\sum_{j=1}^{d} \frac{f_i\left(x+\nu e_{j}\right)-f_i\left(x\right)}{\nu} e_{j}
\end{equation}
where $\nu>0$. From Assumption \ref{ass: smoothness}, we have the following result about the estimation error
\begin{equation}
\Delta^\nu_i(x):={\nabla}^\nu f_i\left(x\right)-\nabla f_i(x).\notag
\end{equation}
\begin{lemma}[\cite{berahas2022theoretical}, Theorem 3.2]
\label{lmm: gradient estimation error}
Under Assumption \ref{ass: smoothness}, we have
\begin{equation}
\left\|\Delta^\nu_i(x)\right\|_{2} \leq \alpha_i \nu, \text{ with } \alpha_i = \frac{\sqrt{d}  M_i}{2}.
\label{eq: gradient approximation error}
\end{equation}
\end{lemma}

\subsection{Local feasible set construction}
\label{sec: Safesetconstruction}
Based on \eqref{eq: gradient approximation} and \eqref{eq: gradient approximation error}, one can build a local feasible set around a strictly feasible point $x_0$ as follows. 
\begin{theorem}
\label{thm: safesetconstruction}
For any strictly feasible point $x_0$, let 
\begin{equation}
\label{eq: l_0}
l_0^* = \min_{i\in\{1, \dots, m\}}\;-f_i(x_0)/L_{\mathrm{max}}, 
\end{equation}
and $\nu_0^* = l_0^*/\sqrt{d}$, where $L_{\mathrm{max}} = \max_{i\geq 1}L_i$. Define 
\begin{align}\label{eq:safe_set}
\mathcal{S}^{(0)}_i&(x_0) :=\big\{x: f_i(x_0)\\\notag
&+{\nabla}^{\nu_0^*} f_i\left(x_0\right)^\top(x-x_0)+2M_i\|x-x_0\|^2\leq 0 \big\}.  
\end{align}
Under Assumption \ref{ass: smoothness}, all the samples needed for computing ${\nabla}^{\nu_0^*} f_i\left(x_0\right)$ are feasible, {$x_0$ belongs to the set $\mathcal{S}^{(0)}(x_0):=\bigcap^m_{i=1}\mathcal{S}^{(0)}_i(x_0)$ and any element of the non-empty $\mathcal{S}^{(0)}(x_0)$ is strictly feasible with respect to \eqref{eq: optimization problem}.}
\end{theorem}

Theorem \ref{thm: safesetconstruction} is similar with \cite[Theorem 1]{guo2022safe}. For completeness, we include the proof of Theorem \ref{thm: safesetconstruction} in Appendix \ref{sec: proofofsafeset}. By construction, we see that if $x_0$ is strictly feasible, then $x_0$ belongs to the interior of $\mathcal{S}^{(0)}(x_0)$ and $\mathcal{S}^{(0)}(x_0)\neq\emptyset$. Moreover, the set $\mathcal{S}^{(0)}(x_0)$ is convex since $\mathcal{S}^{(0)}(x_0) = \cap^{m}_{i=1}\mathcal{S}^{(0)}_i(x_0)$ and $\mathcal{S}^{(0)}_i(x_0)$ is a ball for any $i$. We call $\mathcal{S}^{(0)}(x_0)$ a \textit{local feasible set around} $x_0$.
\begin{remark}
\label{rmk: LM}
The feasibility of $\mathcal{S}^{(0)}(x_0)$ is a consequence of Assumption \ref{ass: smoothness}. Next, we comment on the missing knowledge of $L_i$ and $M_i$ verifying \eqref{eq: LANDM}. In this case, the set $\mathcal{S}^{(0)}(x_0)$ built based on the initial guesses, $L_i$ and $M_i$, might not be feasible. When infeasible samples are generated, one can multiply $L_i$ and $M_i$ for $i\in\mathbb{Z}^m_1$ by  $\beta>1$. This way, at most  $m+\sum^m_{i=1}\max\{\log_\beta(L_{i,\mathrm{inf}}/L_i),\log_\beta(M_{i,\mathrm{inf}}/M_i)\}$ infeasible samples are encountered, where $L_i$ and $M_i$ are the initial guesses. At the same time, one should avoid using a too large value for $M_i$, since if $M_i\gg M_{i,\mathrm{inf}}$, the approximation used to construct $\mathcal{S}^{(0)}(x_0)$ can be very conservative. We refer the readers to Theorem \ref{thm: complexity}, for a discussion on the growth of the complexity of the proposed method with $L_{\max}+M_{\max}$, and Section \ref{sec: experiments}, for an example illustrating the impact of $M_{\max}$ on the convergence.
\end{remark}
{
\begin{remark}
Besides the finite difference method in \eqref{eq: gradient approximation}, there are many other randomized methods for estimating gradient, including Gaussain smoothing \cite{berahas2022theoretical} and Gaussian process regression \cite{shu2022zeroth}. The works \cite{usmanova2022log,berahas2022theoretical,shu2022zeroth} discuss how to decrease the gradient estimation error below a given threshold for a high probability through sampling even when independent and identically distributed measurement noise is present. Note that one can replace ${\nabla}^{\nu_0^*} f_i\left(x_0\right)$ in \eqref{eq:safe_set} with $\widetilde{\nabla}f_i(x_0)$, the gradient estimation obtained using any method. According to the proof of Theorem \ref{thm: safesetconstruction}, we can ensure the feasibility of the local feasible set $\mathcal{S}^{(0)}(x_0)$ for a high probability as long as the estimation error norm verifies $\|\widetilde{\nabla}f_i(x_0) - \nabla f_i(x_0)\| \leq \frac{M_i l^*_0}{2}$ for $i\in \mathbb{Z}^m_1$ for a high probability.
\end{remark}}
One can find in \cite{usmanova2022log} and  \cite{sui2015safe} a different formulation of local feasible sets. In Appendix \ref{sec: comparison of local safe sets}, we compare the two formulations and explain why $\mathcal{S}^{(0)}(x_0)$ is the less conservative.

\subsection{The proposed algorithm}
\label{sec: algorithm}
The proposed method to solve problem \eqref{eq: optimization problem}, called Safe Zeroth-Order Sequential QCQP (SZO-QQ), is summarized in Algorithm~\ref{alg:sslo}. 
The idea is to start from a strictly feasible initial point $x_0$ and iteratively solve \eqref{eq: starstar} in Algorithm~\ref{alg:sslo} until two termination conditions are satisfied. Below, we expand on the main steps of the algorithm.
\begin{algorithm}[htbp!]
\caption{Safe Zeroth-Order Sequential QCQP (SZO-QQ)}
\textbf{Input:} $\mu,\xi,\Lambda>0$, initial feasible point $x_0\in\Omega$\\
\textbf{Output:} $\tilde{x}$, $\tilde{\lambda}$, $\tilde{k}$
\begin{algorithmic}[1]
\State Choose $M_i >M_{i,\mathrm{inf}}$, for $i\in \mathbb{Z}^m_1$
\State $k\gets 0, \textsc{TER}= 0$ 
\While{$\textsc{TER} = 0$}
\State Compute $\mathcal{S}^{(k)}(x_{k})$ based on \eqref{eq:safe_set} and \eqref{eq: stepsize_2}.
\State 
Compute the optimal primal and dual solutions $(x_{k+1},\lambda^\circ_{k+1})$ of 
\begin{equation}
\label{eq: starstar}
\min_{x\in \mathcal{S}^{(k)}(x_{k})} f_0(x)+\mu \|x-x_k\|^2  \tag{SP1}
\end{equation}
\If{$\|x_{k+1}-x_{k}\|\leq \xi$}
\State 
\begin{equation}
\label{eq: SP2}
\lambda_{k+1}\gets \underset{\lambda_{k+1}\in \mathbb{R}^m_{+} 
}{\text{argmin}}   \|\lambda_{k+1}\|_{\infty} \;\;\text{ s.t. }  \eqref{eq: approximateKKT_subproblem} \tag{SP2}
\end{equation}	
\If{$\|\lambda_{k+1}  \|_\infty\leq 2\Lambda$}
\State 
\hspace{-5mm}$\tilde{x}\gets x_{k+1},\tilde{\lambda}\gets \lambda_{k+1}$, $\tilde{k}\gets k+1$, $\textsc{TER}\gets 1$
\EndIf
\EndIf
\State $k\gets k+1$
\EndWhile
\end{algorithmic}
\label{alg:sslo}
\end{algorithm}

\textbf{Providing input data}.
The input to the Algorithm 
\ref{alg:sslo} includes an initial feasible point $x_0$ (see Assumption \ref{ass: strict feasible}) and three parameters $\mu,\xi,\Lambda$. We will describe in Section \ref{sec: accumulation_point} the selection of $\xi$ and $\Lambda$ to ensure that Algorithm \ref{alg:sslo} returns an $\eta$-KKT pair of \eqref{eq: optimization problem}. The impact of $\mu$ on the convergence will be shown in Theorem \ref{thm: complexity}.

\textbf{Building local feasible sets} (Line 4). For a strictly feasible $x_k$, we use \eqref{eq: l_0} to define $l_k^*$ and let the step size of the finite differences for gradient estimation be
\begin{equation}
\nu_k^* =\min\{ \frac{l_k^*}{\sqrt{d}},\frac{1}{k},\frac{\eta}{12\alpha_{\max}m\Lambda}\}.
\label{eq: stepsize_2}
\end{equation} 
Moreover, we use \eqref{eq:safe_set} to define
$\mathcal{S}^{(k)}(x_k)$ in \eqref{eq: starstar}. {The bound $\nu^*_k\leq\frac{l_k^*}{\sqrt{d}}$ ensures that the samples used in gradient estimation are feasible (see Theorem \ref{thm: safesetconstruction}) while $\nu^*_k\leq 1/k$  and  $\nu^*_k\leq{\eta}/(12\alpha_{\max}m\Lambda)$ in \eqref{eq: stepsize_2} are utilized to verify the approximate KKT conditions \eqref{eq:kkt} (see Theorem \ref{prop: eta-KKT}). }

\textbf{Solving a subproblem} (Line 5). Based on the local feasible set, we formulate the subproblem \eqref{eq: starstar} of Algorithm \ref{alg:sslo}. The regulation term $\mu\|x-x_k\|^2$ along with $2M_i\|x-x_k\|$ in $\mathcal{S}^{(k)}_i(x_k)$ prevents too large step sizes. If $\|x-x_k\|$ is large, the proxies used in \eqref{eq: starstar} are not accurate. The regulation term can also be found in the proximal trust-region method in \cite{aravkin2022proximal}. With it, we can ensure that $\|x_{k+1}-x_k\|$ converges to 0 (see Proposition \ref{thm: main_theorem}) and conduct complexity analysis (see Theorem \ref{thm: complexity}). Since $f_0$ is assumed, without loss of generality, to be known and linear (see Section \ref{sec::Prob}), \eqref{eq: starstar} is a known convex QCQP. {As an extension to Theorem \ref{thm: safesetconstruction}, we have the following proposition, the proof of which is in Appendix \ref{prf: feasibility of sp}}.
{
\begin{proposition}
\label{prop:feasibility of (SP)}
Under Assumptions \ref{ass: smoothness} and \ref{ass: strict feasible}, for any $k$, $x_k$ is a strictly feasible solution to $\mathrm{(SP1)}$ and to \eqref{eq: optimization problem}.
\end{proposition}
}
Although for any $k$ the iterate $x_k$ is strictly feasible w.r.t. \eqref{eq: optimization problem}, it is possible that a subsequence of the sequence $\{x_k\}_{k\geq 1}$ converges to a point on the boundary of $\Omega$.

\textbf{Checking termination conditions} (Line 6-11). We introduce two termination conditions guaranteeing that the pair $(\tilde{x},\tilde{\lambda})$ returned by Algorithm \ref{alg:sslo} is an $\eta$-KKT pair. The first one (Line 6) requires that $\|x_{k+1}-x_k\|$ is smaller than a given threshold $\xi$ while the second requires that the solution to the optimization problem \eqref{eq: SP2} is small enough (Line 8). The constraint of \eqref{eq: SP2} is 
\begin{equation}
\label{eq: approximateKKT_subproblem}
\max\Big\{\delta_1(k,\lambda_{k+1}), \max\{\delta^{(i)}_2(k,\lambda_{k+1}):i\geq 1 \}\Big\}\leq\frac{\eta}{2},
\end{equation}
where
\begin{align}\notag
\delta_1(k,\lambda_{k+1})& :=  \Big\|{\nabla} f_0\left(x_{k+1}\right) + 2\mu (x_{k+1}-x_k) \\\notag
 +\sum^{m}_{i=1} \lambda_{k+1}^{(i)}&\left({\nabla}^{\nu_{k}^*} f_i\left(x_{k}\right)+4M_i(x_{k+1}-x_{k})\right)\Big\|,\\\notag
\delta^{(i)}_2(k,\lambda_{k+1}) &:=\Big|\lambda_{k+1}^{(i)} \Big(f_i(x_k)+{\nabla}^{\nu_k^*} f_i\left(x_k\right)(x_{k+1}-x_k)\\\label{eq: definition_delta}
&\qquad\qquad+2M_i\|x_{k+1}-x_k\|^2\Big)\Big|.
\end{align}
Observe that $\delta_1(k,\lambda_{k+1})$ and $\delta^{(i)}_2(k,\lambda_{k+1})$ in \eqref{eq: definition_delta} originate from the KKT conditions for \eqref{eq: starstar}. Therefore, by solving \eqref{eq: SP2} we obtain the smallest-norm vector $\lambda_{k+1}$ such that $(x_{k+1},\lambda_{k+1})$ is a $\eta/2$-KKT pair of \eqref{eq: starstar}. To solve \eqref{eq: SP2}, we reformulate it as a convex QCQP and use $\lambda^\circ_{k+1}$ (see Line 5 of Algorithm \ref{alg:sslo}) as an initial feasible solution. If the two conditions are satisfied at the $(k+1)$-th iteration, then the algorithm outputs in Line 9 are $\tilde{x} = x_{k+1}$, $\tilde{\lambda} = \lambda_{k+1}$ and $\tilde{k} = k+1$.

Algorithm \ref{alg:sslo} is similar to Sequential QCQP (SQCQP)~\cite{fukushima2003sequential}. In SQCQP, at each iteration, quadratic proxies for constraint functions are built based on the local gradient vectors and Hessian matrices. The application of SQCQP to optimal control has received increasing attention \cite{messerer2020determining,messerer2021survey}, due to the development of efficient solvers for QCQP subproblems \cite{frison2022introducing}. Different from SQCQP \cite{fukushima2003sequential}, Algorithm~\ref{alg:sslo} can guarantee sample feasibility and does not require the knowledge of Hessian matrices, which are costly to obtain for zeroth-order methods. As Hessian matrices are essential for proving the convergence of SQCQP in \cite{fukushima2003sequential}, we cannot use the same arguments in \cite{fukushima2003sequential} to show the properties of SZO-QQ's iterates. In the following two sections, we state the properties of $(\tilde{x},\tilde{\lambda})$ and analyze the efficiency of the algorithm.  

\section{Properties of SZO-QQ's iterates and output} 
\label{sec: accumulation_point}
In this section, we aim to show that, for a suitable $\xi$, the pair $(\tilde{x},\tilde{\lambda})$ derived in Algorithm~\ref{alg:sslo} is $\eta$-KKT. We start by considering the infinite sequence of Algorithm~\ref{alg:sslo}'s iterates $\{x_k\}_{k\geq 1}$ when the termination conditions in Line 6 and 8 of Algorithm \ref{alg:sslo} are removed. We show that the sequence $\{x_k\}_{k\geq 1}$ has accumulation points and, for any accumulation point $x_c$, under mild assumptions, there exists $\lambda_c\in\mathbb{R}^m_{\geq 0}$ such that $(x_c,\lambda_c)$ is a KKT pair of \eqref{eq: optimization problem}. Based on this result, we then study the activation of the two termination conditions and prove that they are satisfied within a finite number of iterations. Finally, we show that if $\xi$ is carefully chosen, the derived pair $(\tilde{x},\tilde{\lambda})$ is $\eta$-KKT. 

\subsection{On the accumulation points of $\{x_k\}_{k\geq 1}$}
\label{sec: convergence of x_K}
\begin{proposition}
\label{thm: main_theorem}
If the termination conditions are removed, the sequence $\{x_k\}_{k\geq 1}$ in Algorithm~\ref{alg:sslo} has the following properties:
\begin{itemize}
\item[1.] the sequence $\{f_0(x_k)\}_{k\geq 1}$ is non-increasing;
\item[2.] $\{x_k\}_{k\geq 1}$ has at least one accumulation point $x_c$ and $\{\|x_{k+1}-x_k\|\}_{k\geq 1}$ converges to 0;
\item [3.] $\lim_{k\rightarrow \infty}f_0(x_k)=f_0(x_c)> -\infty$. 
\end{itemize}
\end{proposition}
The proof is provided in Appendix \ref{sec: 1st argument of thm2}. It mainly exploits the following inequality,
\begin{equation}
\label{eq: monotonicity_inequality}
f_0(x_{k+1})+\mu\|x_{k+1}-x_k\|^2\leq f_0(x_{k}),
\end{equation}
which is due to the optimality of $x_{k+1}$ for~\eqref{eq: starstar} in Algorithm \ref{alg:sslo}. The monotonicity of $f_0(x_k)$ and the convergence of $\|x_{k+1}-x_k\|$ are direct consequences of \eqref{eq: monotonicity_inequality}. By utilizing the monotonicity, we have that, for any $k\geq 1$, $x_k$ belongs to the bounded set $ \mathcal{P}_\beta$ (see Assumption \ref{ass: bounded_sublevel} for the definition of $\mathcal{P}_\beta$). Due to Bolzano–Weierstrass theorem, there exists an accumulation point of $\{x_k\}_{k\geq 1}$. The continuity of $f_0(x)$ gives us Point 3 of Proposition \ref{thm: main_theorem}. 

Based on Proposition \ref{thm: main_theorem}, we can show that under Assumption \ref{ass:LICQ} below, there exists an accumulation point of $\{x_k\}_{k\geq 1}$ that allows one to build a KKT pair.
\begin{assumption}
\label{ass:LICQ}
There exists an accumulation point $x_c$ of $\{x_k\}_{k\geq 1}$ such that 
the~\textit{Linear Independent Constraint Qualification (LICQ)} holds at $x_c$ for \eqref{eq: optimization problem}, which is to say the gradients $\nabla f_i(x_c)$ with $i\in \mathcal{A}(x_c):=\{i: f_i(x_c)=0\}$ are linearly independent.
\end{assumption}
Assumption \ref{ass:LICQ} is widely used in optimization \cite{wachsmuth2013licq}. For example, it is used to prove the properties of the limit point of the Interior Point Method \cite{jorge2006numerical}. With this assumption, if $x_c$ is a local minimizer, then there exists $\lambda_c\in\mathbb{R}^{m}_{\geq 0}$ such that $(x_c,\lambda_c)$ is a KKT pair \cite[Theorem 12.1]{jorge2006numerical}, which will be used in the proof of Theorem \ref{thm: accumulationKKT}.
\begin{theorem}
\label{thm: accumulationKKT}
Let Assumption \ref{ass:LICQ} hold, and let $x_c$ be an accumulation point of $\{x_k\}_{k\geq 1}$ where LICQ is verified. Then, there exists a unique $\lambda_c\in \mathbb{R}^m_{\geq0}$ such that $(x_c,\lambda_c)$ is a KKT pair of the problem \eqref{eq: optimization problem}.
\end{theorem}
We only consider the case where $\mathcal{A}(x_c)$ is not empty in the proof of Theorem \ref{thm: accumulationKKT}, provided in Appendix \ref{sec: proof_of_accumulationKKT}. The proof can be easily extended to account for the case where $x_c$ is in the interior of $\Omega$. To show Theorem \ref{thm: accumulationKKT}, we exploit a preliminary result (Lemma \ref{lmm: non_empty_interior}, stated and proved in Appendix \ref{sec: preliminary results for Lemma xcandlambdac}) where we construct an auxiliary problem \eqref{eq: x_c_optimization} and show that $x_c$ is an optimizer to \eqref{eq: x_c_optimization}. We notice that the KKT conditions of \eqref{eq: x_c_optimization} evaluated at $x_c$ coincide with those of \eqref{eq: optimization problem} evaluated at the same point. Due to LICQ, there exists a unique $\lambda_c\in \mathbb{R}_{\geq0}$ such that $(x_c,\lambda_c)$ is a common KKT pair of \eqref{eq: x_c_optimization} and \eqref{eq: optimization problem} \cite[Section 12.3]{jorge2006numerical}.
{
\begin{remark}
Theorem \ref{thm: accumulationKKT} extends the asymptotic results in \cite{guo2022safe} to solve the problems with non-convex objective and constraint functions. Our proof for Theorem 2 utilizes the marginal difference between the KKT conditions of the original problem and the QCQP subproblems. This approach is different from the one used in \cite{guo2022safe}, where the authors used convexity of the feasible set to prove the feasibility of a line segment by showing that the two endpoints are feasible. Worth noticing, this extension does not come for free either considering that Assumption 4 in the current paper is slightly more stringent than Assumption 4 in \cite{guo2022safe}.
\end{remark}}
\subsection{The output of Algorithm \ref{alg:sslo} is an $\eta$-KKT pair}
\label{sec: eta-KKT}
The result in Theorem \ref{thm: accumulationKKT} is asymptotic, but in practice, only finitely many iterations can be computed. From now on, we take the termination conditions of Algorithm~\ref{alg:sslo} into account and show that, given any $\eta>0$, by suitably tuning $\xi>0$, Algorithm \ref{alg:sslo} returns an $\eta$-KKT pair. First, we make the following assumption, which allows us to conclude in Proposition \ref{thm: algorithm_termination} the finite termination of Algorithm \ref{alg:sslo}.
\begin{assumption}
\label{ass: magnitude_dual_variable}
The KKT pair $(x_c,\lambda_c)$ in Theorem~\ref{thm: accumulationKKT} satisfies $\|\lambda_c\|_{\infty}<\Lambda$, where $\Lambda >0$ is the input in Algorithm \ref{alg:sslo}.
\end{assumption}
Assumptions on the bound of the dual variable are adopted in the literature of primal-dual methods, including~\cite{usmanova2021fast,zhu2016distributed}. We illustrate Assumption \ref{ass: magnitude_dual_variable} in Appendix~\ref{sec: characterization of lambda} where we show in an example that $\Lambda$ is related to the geometric properties of the feasible region. {In case it is hard to find $\Lambda$ satisfying Assumption \ref{ass: magnitude_dual_variable}, we refer to Remark \ref{rmk: on Lambda} on updating the initial guess of $\Lambda$.}

\begin{proposition}
\label{thm: algorithm_termination}
Let Assumptions~\ref{ass:LICQ} and~\ref{ass: magnitude_dual_variable} hold, Algorithm~\ref{alg:sslo} terminates in a finite number of iterations.
\end{proposition}
According to Proposition~\ref{thm: main_theorem}, the first termination condition is satisfied in Algorithm~\ref{alg:sslo} whenever $k$ is sufficiently large. In the proof of Proposition~\ref{thm: algorithm_termination} (provided in Appendix \ref{sec: algorithm_termination}), we show that $\lambda_{k+1}=\lambda_c$ is a feasible solution to \eqref{eq: SP2} when $x_{k+1}$ is close enough to $x_c$. Thus, for sufficiently large $k$, the second termination is satisfied since $\|\lambda_c\|_\infty<\Lambda$ .

Recall that Algorithm \ref{alg:sslo} returns $\tilde{x},\tilde{\lambda}$ and $\tilde{k}$, which are dependent on the chosen value for $\xi$. For a given accuracy indicator $\eta>0$, in the following, we show how to select $\xi$ such that $(\tilde{x},\tilde{\lambda})$ is an $\eta$-KKT pair.
\begin{theorem}
\label{prop: eta-KKT}
Let Assumption \ref{ass:LICQ} hold, and let
\begin{align}\label{eq: xi}
\xi=h(\eta) :=\min\Big\{ &\frac{\eta}{60\Lambda\sum^m_{i=1}M_i},\frac{\eta}{12\mu}, 1,\\ \notag
&\frac{\eta}{4\Lambda(   \alpha_{\max}+2{L_{\max}} + 2{M_{\max}})}\Big\}
\end{align}
be satisfied, where $\mu$ is a parameter of \eqref{eq: starstar}, $\alpha_{\max} = \max_{1\leq i\leq m} \alpha_i$ and $\alpha_i$ is defined in \eqref{eq: gradient approximation error}.
Then the output $(\tilde{x},\tilde{\lambda})$ of Algorithm \ref{alg:sslo} is an $\eta$-KKT pair of \eqref{eq: optimization problem} .
\end{theorem}
The proof of Theorem \ref{prop: eta-KKT} can be found in Appendix \ref{sec: proof_of_xi}. 

\begin{remark}
\label{rmk: on Lambda}
In case it is hard to choose a value of $\Lambda$ fulfilling Assumption \ref{ass: magnitude_dual_variable}, we can replace $\Lambda$ with $\kappa\|\lambda_{k+1}\|_\infty$, where $\kappa>1$ and $\lambda_{k+1}$ is the solution to the problem \eqref{eq: SP2} in Algorithm \ref{alg:sslo}, every time the second termination condition (Line 8 in Algorithm \ref{alg:sslo}) is violated. Note that every time $\Lambda$ gets updated, it becomes at least $2\kappa-1$ times larger. Similar updating rules can also be found in~\cite{usmanova2021fast}. {In this way, Algorithm \ref{alg:sslo} is bound to terminate in a finite number of iterations. Otherwise, the second termination condition would be violated for infinitely many times, $\Lambda$ would keep growing and we would find $\Lambda$ that satisfies Assumption~\ref{ass: magnitude_dual_variable}. With Assumption~\ref{ass: magnitude_dual_variable} satisfied, Algorithm \ref{alg:sslo} would terminate in a finite number of iterations (as indicated in Proposition \ref{thm: algorithm_termination}). }

{Note that when Algorithm \ref{alg:sslo} terminates, the updated $\Lambda$ does not necessarily satisfy Assumption \ref{ass: magnitude_dual_variable}. However, one can check the proof and verify that Theorem \ref{prop: eta-KKT} does not rely on Assumption \ref{ass: magnitude_dual_variable}. Therefore, even if we fail to obtain $\Lambda$ that satisfies Assumption \ref{ass: magnitude_dual_variable}, Algorithm \ref{alg:sslo} returns an $\eta$-KKT pair if $\xi$ is selected according to \eqref{eq: xi}.  }

However, we also notice that this updating mechanism might generate a conservative guess for $\Lambda$ if $\|\lambda_{k}\|_\infty\gg \|\lambda_c\|_\infty$ for some $k$. In Theorem \ref{prop: eta-KKT}, we will set $\xi$ in Algorithm \ref{alg:sslo} to be proportional to $\Lambda^{-1}$ so that the returned pair is an $\eta$-KKT pair. Consequently, a conservative $\Lambda$ can increase the number of iterations required by Algorithm \ref{alg:sslo}.
\end{remark}

In summary, to ensure that the pair $(\tilde{x},\tilde{\lambda})$ is an $\eta$-KKT pair, we need to set $\xi$ in {Algorithm}~\ref{alg:sslo} to be $h(\eta)$ in \eqref{eq: xi} while selecting $L_i$, $M_i$ and $\Lambda$ for satisfaction of Assumptions \ref{eq: smoothness} and \ref{ass: magnitude_dual_variable} or according to Remarks \ref{rmk: LM} and \ref{rmk: on Lambda}.

\section{Complexity analysis}
\label{sec: complexity_analysis}
In this section, we aim to give an upper bound, dependent on $\eta$, for the number of Algorithm~\ref{alg:sslo}'s iterates. To this purpose, we consider the following assumption, which allows us to show the convergence of $\{x_k\}_{k\geq 1}$ in Lemma~\ref{lmm: convergence}.
\begin{assumption}
\label{ass: assumption for convergence}
The accumulation point $x_c$ in Assumption \ref{ass:LICQ}, which is already known to be the primal of some KKT pair, is a strict local minimizer, i.e., there exists a neighborhood $\mathcal{N}$ of $x_c$ such that $f_0(x)>f_0(x_c)$ for any $x\in \mathcal{N}\cap \Omega\setminus x_c$.
\end{assumption}
If Second-Order Sufficient Condition (SOSC) for optimality \cite[Chapter 12.5]{jorge2006numerical} is satisfied, Assumption \ref{ass: assumption for convergence} holds. Since this assumption does not rely on the twice differentiability of the objective and constraint functions, it is more general than SOSC, which is commonly assumed in the optimization literature~\cite{diehl2005real,izmailov2004newton}. 
\begin{lemma}
\label{lmm: convergence}
If Assumption \ref{ass: assumption for convergence} holds, $\{x_k\}_{k\geq 1}$ converges.
\end{lemma}
The proof of Lemma \ref{lmm: convergence} is in Appendix \ref{sec: proof of convergence lemma}. In the rest of this section, we consider Assumption \ref{ass: assumption for convergence}. Let $x_c$ be the limit point of $\{x_k\}_{k\geq 1}$ and note that there exists $\lambda_c$ such that $(x_c,\lambda_c)$ is a KKT pair. Then we show in Lemma \ref{lmm: convergence of the dual variable}, with the proof in Appendix 
\ref{sec: lambda_conv}, that $\lambda^\circ_{k}$, the optimal dual solution to \eqref{eq: starstar}, converges to $\lambda_c$. 
\begin{lemma}
\label{lmm: convergence of the dual variable}
Let Assumptions \ref{ass:LICQ}, \ref{ass: magnitude_dual_variable} and \ref{ass: assumption for convergence} hold, $\lambda^\circ_{k}$ converges to $\lambda_c$. 
\end{lemma}
With Lemma \ref{lmm: convergence of the dual variable} and Assumption~\ref{ass: magnitude_dual_variable}, we know that for any $\eta$ there exists $K$, independent of $\eta$, such that $\|\lambda_k\|_\infty\leq 2\Lambda$ for any $k\geq K$. Therefore, for sufficiently small $\eta$, Algorithm \ref{alg:sslo} terminates whenever the first termination condition in Line 6 is satisfied. We can now conclude in Theorem \ref{thm: complexity} on the complexity of Algorithm \ref{alg:sslo} by analyzing only the first termination condition. The proof of Theorem~\ref{thm: complexity} is in Appendix~\ref{sec: proof_of_complexity}.

\begin{theorem}
\label{thm: complexity}
Let Assumptions \ref{ass:LICQ}, \ref{ass: magnitude_dual_variable} and \ref{ass: assumption for convergence} hold, there exists {$0<\bar{\eta}<1$} such that, for any $\eta<\bar{\eta}$, Algorithm \ref{alg:sslo} terminates within $\overline{\mathcal{K}}(\eta)+1$ iterations, where 
\begin{equation}
\label{eq: iteration_upperbound}
\overline{\mathcal{K}}(\eta) = \frac{f_0(x_0) - \inf\{f_0(x): x\in \Omega\}}{\mu\cdot (h(\eta))^2},    
\end{equation}
and $\mu$ is the coefficient of the quadratic penalty term in \eqref{eq: starstar}, and $h(\eta)$ is defined in \eqref{eq: xi}. {Thus, given $C_1>0$, if the initial point is properly selected such that $f_0(x_0) - \inf\{f_0(x): x\in \Omega\}<C_1$, then for any {$0<\eta\leq \bar{\eta}$}, Algorithm~\ref{alg:sslo} takes at most  $O\left({{d}}/{\eta^2}\right)$
iterations and $O\left({d^{2}}/{\eta^2}\right)$ samples to return $(\tilde{x},\tilde{\lambda})$, an $\eta$-KKT pair of the problem~\eqref{eq: optimization problem}}.
\end{theorem}
{We notice that $\overline{\mathcal{K}}(\eta)$ is a non-smooth univariate function of the hyperparameter $\mu$ of Algorithm \ref{alg:sslo}. If $\mu$ is too large, according to \eqref{eq: xi}, $\xi = h(\eta)$ in the termination condition (Line 6 of Algorithm \ref{alg:sslo}) can be too small, resulting in unacceptably large $\overline{\mathcal{K}}(\eta)$. If $\mu$ is too small, $h(\eta)$ in \eqref{eq: xi} is independent of $\mu$ and thus $\overline{\mathcal{K}}(\eta)$, inversely proportional to $\mu$, can again be very large. Therefore, there exists at least one local minimizer for $\overline{\mathcal{K}}(\eta)$. For finding it, one can use pattern search \cite{audet2002analysis}.}  

\textbf{Discussion:}
We compare the sample and computation complexities of SZO-QQ with two other existing safe zeroth-order methods, namely, LB-SGD in \cite{usmanova2022log} and SafeOpt in \cite{berkenkamp2016safe}. We remind the readers that these methods have different assumptions. Specifically, given the black-box optimization problem \eqref{eq: optimization problem}, SafeOpt assumes that $f_i(x), i\in \mathbb{Z}^m_{i=0}$, has bounded norm in a suitable Reproducing Kernel Hilbert Space while both SZO-QQ and LB-SGD assume the knowledge of the smoothness constants. {Regarding sample complexity, SZO-QQ needs $O(\frac{d^2}{\eta^2})$ samples under mild conditions to generate an $\eta$-KKT pair while LB-SGD and SafeOpt require at least 
$O(\frac{d^2}{\eta^5})$  and $O(\frac{\rho(d)}{\eta^2})$ samples\footnote{In \cite{usmanova2022log}, the authors claim the iterations of LB-SGD needed for obtaining an $\eta$-KKT to be $O(\frac{d^2}{\eta^5})$ when the objective and constraint functions are strongly convex. In \cite[Theorem 1]{sui2015safe}, $O(\frac{\rho(d)}{\eta^2})$ samples are needed in SafeOpt to obtain an $\eta$-suboptimal point. For unconstrained optimization with a strongly convex objective function, there exists $\alpha>0$ such that any $\eta$-KKT pair has a $\alpha\eta$-suboptimal primal. Therefore, SafeOpt requires at least $O(\frac{\rho(d)}{\eta^2})$ samples to obtain an $\eta$-KKT pair of \eqref{eq: optimization problem}.  } respectively, where $\rho(d)$ relates to the discretization of the space and therefore can be exponential with $d$.} 

We also highlight that the computational complexity of each iteration of both LB-SGD and SZO-QQ stays fixed while the computation time required for the Gaussian Process regression involved in each iteration of SafeOpt increases as the data set gets larger. The high computational cost is one of the main reasons why SafeOpt scales poorly to high-dimensional problems. Numerical results comparing the computation time and the number of samples required by these methods are provided in Section~\ref{sec: toy example}.

In contrast, the Interior Point Method, based on the assumption of $f_i(x)$ being twice continuously differentiable, achieves superlinear convergence \cite{jorge2006numerical} by utilizing the true model of the optimization problem, which translates into at most $O(\mathrm{log}\frac{1}{\eta})$ iterations. The gap between $O(\frac{1}{\eta^2})$ of SZO-QQ and $O(\mathrm{log}\frac{1}{\eta})$ may be either the price we pay for the lack of the first-order and second-order information of the objective and constraint functions or due to the conservative analysis in Theorem \ref{thm: complexity}. To see whether there exists a tighter complexity bound than $O(\frac{1}{\eta^2})$ for Algorithm \ref{alg:sslo}, an analysis on the convergence rate is needed, which is left as future work.

\section{Numerical Results}
\label{sec: experiments}
In this section, we present three numerical experiments to test the performance of Algorithm \ref{alg:sslo}. The first is a two-dimensional problem where we compare SZO-QQ with other existing zeroth-order methods and discuss the impact of parameters. In the remaining two examples, we apply our method to solve optimal control and optimal power flow problems, which have more dimensions and constraints. All the numerical experiments have been executed on a PC with an Intel Core i9 processor.  {For solving (SP1) and (SP2) in Algorithm \ref{alg:sslo} (the latter can be reformulated as a QCQP), we use
Gurobi~\cite{gurobi}, called through MATLAB via YALMIP~\cite{Lofberg2004}. Since the subproblems of different iterations share the same structure, we use the function \texttt{optimizer} of Yalmip to enhance the interfacing efficiency. A matlab implementation of the following experiments is available at https://github.com/odetojsmith/SZO-QQ.}

\subsection{Solving an unknown non-convex QCQP}
\label{sec: toy example}
We evaluate SZO-QQ and compare it with alternative safe zeroth-order methods in the following non-convex example,
\begin{equation}
\label{eq: toy_example}    
    \begin{aligned}
    \min_{x \in \mathbb{R}^2} & \quad  f_0(x)=0.1\times (x^{(1)})^2+x^{(2)}\\
     \text{subject to}  & \quad f_1(x) =0.5-\|x+[0.5 \; -0.5]^\top\|^2 \leq 0,\\
     &\quad f_2(x)=x^{(2)}-1\leq 0,\\
    & \quad f_3(x) = (x^{(1)})^2-x^{(2)}\leq 0.
    \end{aligned}
\end{equation}
We assume that the functions $f_i(x), i=1,2,3,$ are unknown but can be queried. A strictly feasible initial point $x_0 = [0.9 \,\,\, 0.9]^\top$ is given. The unique optimum  $x_* = [0\,\,\,0]^\top $ is not strictly feasible. According to Theorem \ref{thm: accumulationKKT}, the iterates of SZO-QQ will get close to $x_*$, which allows us to see whether SZO-QQ stays safe and whether the convergence is fast when the iterates are close to the feasible region boundary. The experiment results allow us to discuss the following three aspects, respectively on the derivation of a $10^{-2}$-KKT pair, performance evaluation, and parameter tuning.

\subsubsection{Selection of $\xi$ for deriving a $10^{-2}$-KKT pair}

To begin, we fix $\eta = 10^{-2}$ and aim to derive an $\eta$-KKT pair. By setting $\Lambda = 1.5$ and $L_i=5$, $M_i=3$ for any $i\geq 1$, we calculate $\xi=1.51\times 10^{-5}$ according to~\eqref{eq: xi}. With these values, SZO-QQ returns in 3.3 seconds an $\eta'$-KKT pair with $\eta'=9.21\times 10^{-4}<\eta$. {We notice that all iterates and samples are feasible points for \eqref{eq: toy_example}. This is due to Proposition \ref{prop:feasibility of (SP)} and the fact that the results obtained by solving (SP1) using Gurobi are feasible points for (SP1).} We also observe that $\|{\lambda}_{k}\|_\infty$ converges to 1 and thus the original guess $\Lambda=1.5$ satisfies Assumption \ref{ass: magnitude_dual_variable}. Now we see that we indeed derive an $\eta$-KKT pair, which coincides with Theorem \ref{prop: eta-KKT}. To further evaluate the performance of SZO-QQ in terms of how fast the objective function value decreases, we eliminate the termination conditions in the remainders of Section \ref{sec: toy example}.

\subsubsection{Performance comparison with other methods}

We run SZO-QQ with $\xi=0$ and compare with LB-SGD \cite{usmanova2022log}, Extremum Seeking \cite{arnold2015model} and SafeOptSwarm\footnote{SafeOptSwarm is a variant of SafeOpt (recall Section \ref{sec: introduction}). The former add heuristics to make SafeOpt in \cite{sui2015safe} more tractable for higher dimensions.}~\cite{berkenkamp2016safe}. Among these methods, SZO-QQ, LB-SGD, and SafeOpt-Swarm have theoretical guarantees for sample feasibility (at least with a high probability). Only SZO-QQ and LB-SGD require Assumption \ref{ass: smoothness} on Lipschitz and smoothness constants. For these two approaches, by trial and error (see Remark \ref{rmk: LM}), we choose $L_{i} = 5$ and $M_{i}= 3$ for any $i\geq 1$. The penalty coefficient $\mu$ of Algorithm \ref{alg:sslo} in \eqref{eq: starstar} is set to be $ 0.001$. For both LB-SGD and Extremum Seeking are barrier-function-based, we use the reformulated unconstrained problem $\min_x f_{\mathrm{log}}(x,\mu_{\log})$, where $f_{\mathrm{log}}(x,\mu_{\log}):=f_0(x) - \mu_{\mathrm{log}}\sum^4_{i=1}\log(-f_i(x))$, and $\mu_{\log} = 0.001$.
\begin{figure}[htbp!]
    \centering
    \includegraphics[width=0.6\linewidth]{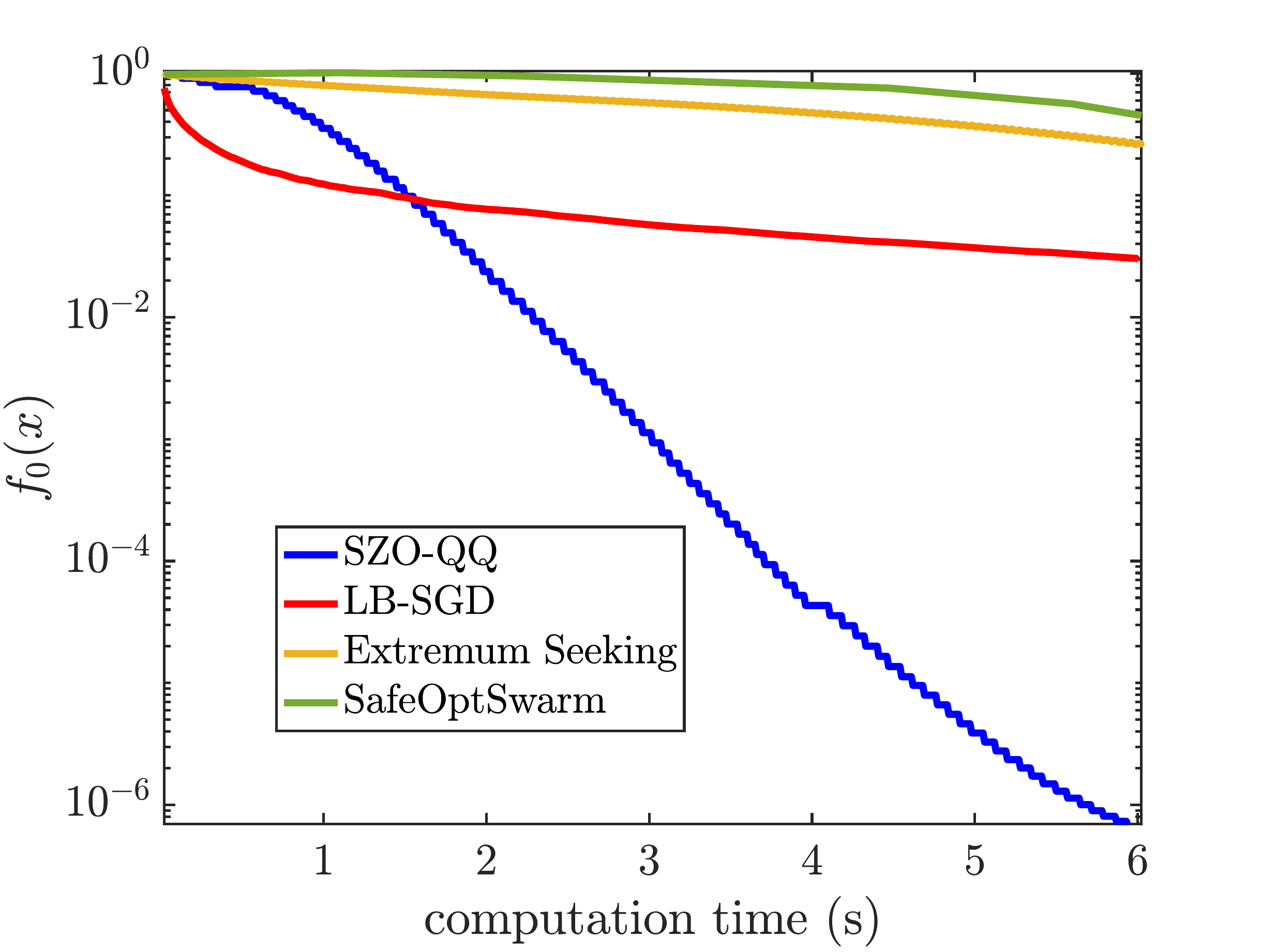}
    \caption{Objective value as a function of computation time.}
    \label{fig:objective_decrease}
\end{figure}

\begin{figure}[htbp!]
    \centering
    \includegraphics[width=0.6\linewidth]{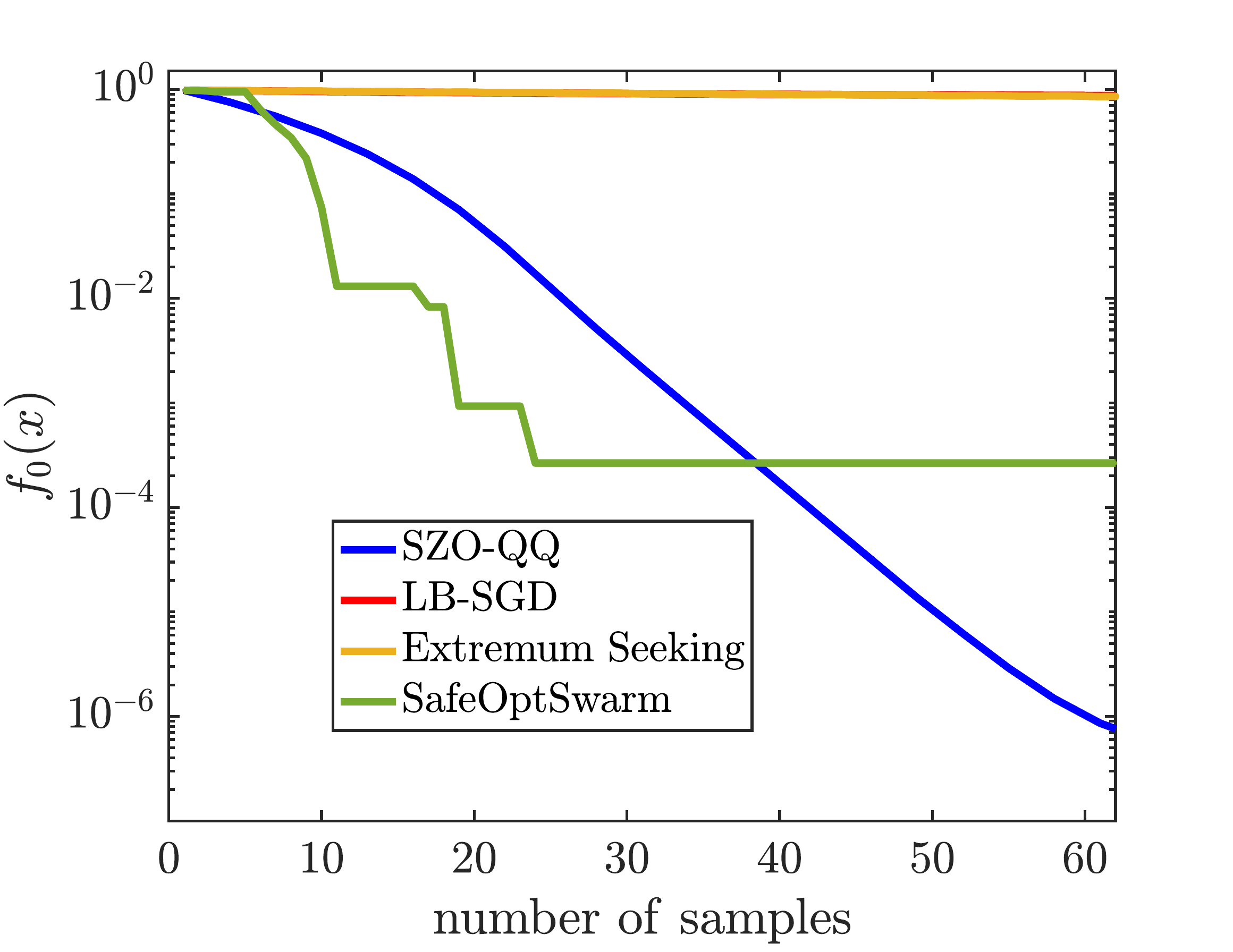}
    \caption{Objective value as a function of the number of samples, where the curve for LB-SGD overlaps with that for Extremum Seeking.}
    \label{fig:objective_decrease_wrt_n_samples}
\end{figure}

In Fig. \ref{fig:objective_decrease}, we show the objective function values versus the computation time {(excluding the time used for sampling)}. During the experiments, none of the methods violates the constraints. Regarding the convergence to the minimum, we see that LB-SGD has the most decrease in the objective function value in the first 1.5 seconds due to the low complexity of each iteration. {In these 1.5 seconds, LB-SGD utilizes 32546 function samples while SZO-QQ only 21.} Afterwards, SZO-QQ achieves a better solution. In the first 6 seconds, SZO-QQ shows a clear convergence trend to the optimum, consistent with Theorem~\ref{thm: accumulationKKT}. {In Fig. \ref{fig:objective_decrease_wrt_n_samples}, we show the objective function values versus the number of samples. We see that SZO-QQ and SafeOptSwarm are much more sample-efficient than LB-SGD and Extremum Seeking. However, to generate the iterates shown in Fig. \ref{fig:objective_decrease_wrt_n_samples}, SafeOptSwarm needs 82 seconds while SZO-QQ only uses 6 seconds.}

LB-SGD slows down when the iterates are close to the boundary of the feasible set (see Appendix \ref{sec: comparison of local safe sets} for the explanation for this phenomenon). Meanwhile, the slow convergence of Extremum Seeking is due to its small learning rate. If the learning rate is large, the iterates might be brought too close to the boundary of the feasible set, and then the perturbation added by this method would lead to constraint violation. These considerations constitute the main dilemma in parameter tuning for Extremum Seeking. Meanwhile, exploring the unknown functions in SafeOptSwarm is based on Gaussian Process (GP) regression models instead of local perturbations. Since SafeOptSwarm does not exploit the convexity of the problem, it maintains a safe set and tries to expand it to find the global minimum. {Empirically, this method samples many points close to the boundary of the feasible region even though they are far from the optimizer, which is also observed in \cite{turchetta2019safe}. This phenomenon is due to the tendency of SafeOptSwarm to enlarge the estimated safe set, which is one of the main reasons why in Fig. \ref{fig:objective_decrease_wrt_n_samples} SafeOptSwarm fails to further decrease the objective function value after 25 samples.} We also run LB-SGD and Extremum Seeking with different penalty coefficients $\mu_{\mathrm{log}}$ to check whether the slow convergence is due to improper parameter tuning. We see that with larger $\mu_{\mathrm{log}}$ the performance of the log-barrier-based methods deteriorates. This is probably because the optimum of the unconstrained problem $f_{\log}(x,\mu_{\log})$ deviates more from the optimum as $\mu_{\mathrm{log}}$ increases. With smaller $\mu_{\mathrm{log}}$,  the Extremum Seeking method leads to constraint violation while the performance of LB-SGD barely changes. 

\subsubsection{Impact of the parameters $L_i$ and $M_i$}

\begin{figure}[htbp!]
    \centering
    \includegraphics[width=0.6\linewidth]{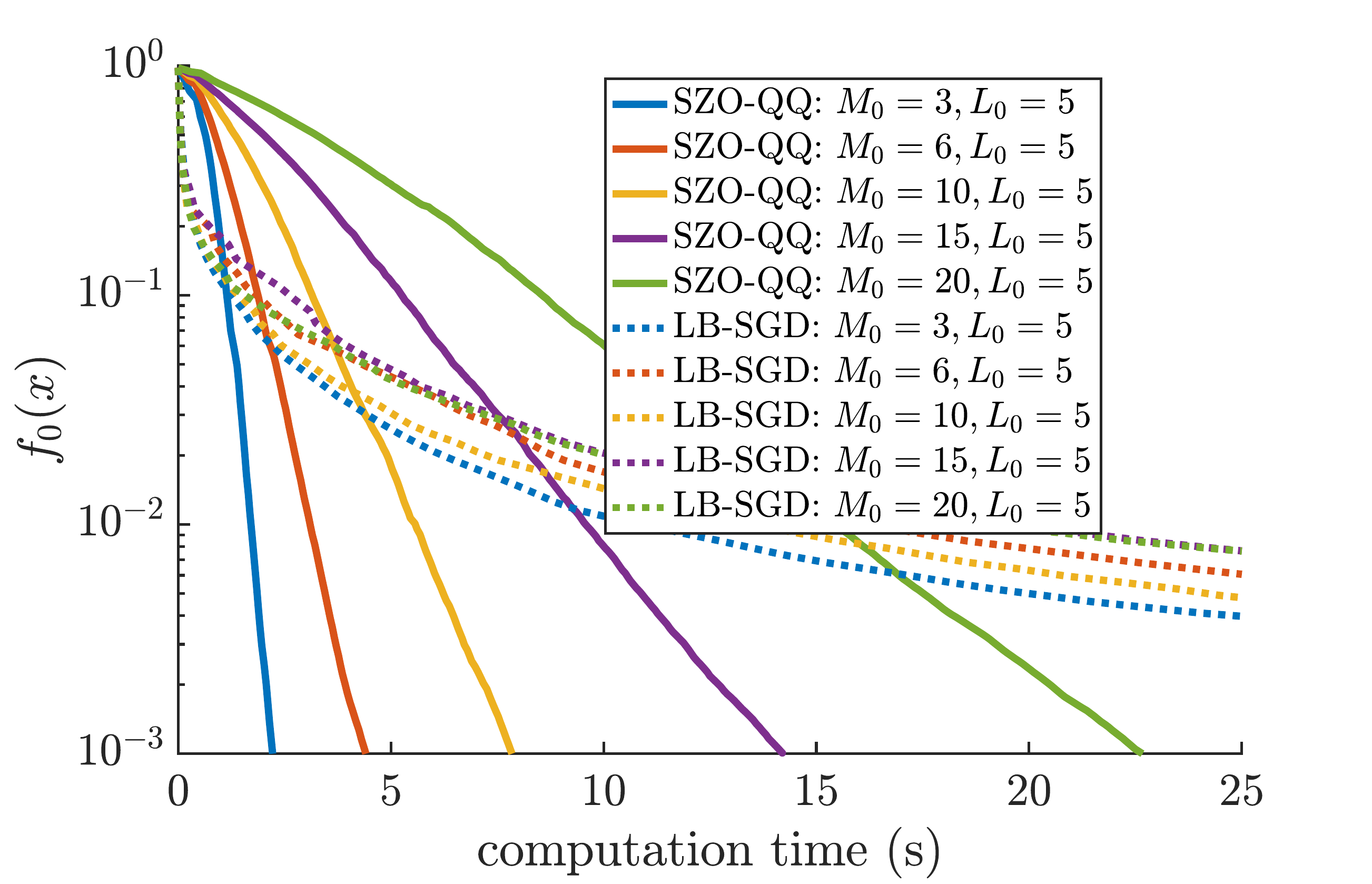}
    \caption{{Objective value as a function of computation time:$\text{ }$ $ L_0=5$ fixed and $M_0$ varied}}
    \label{fig:objective_decrease:M fixed}
\end{figure}

\begin{figure}[htbp!]
    \centering
    \includegraphics[width=0.6\linewidth]{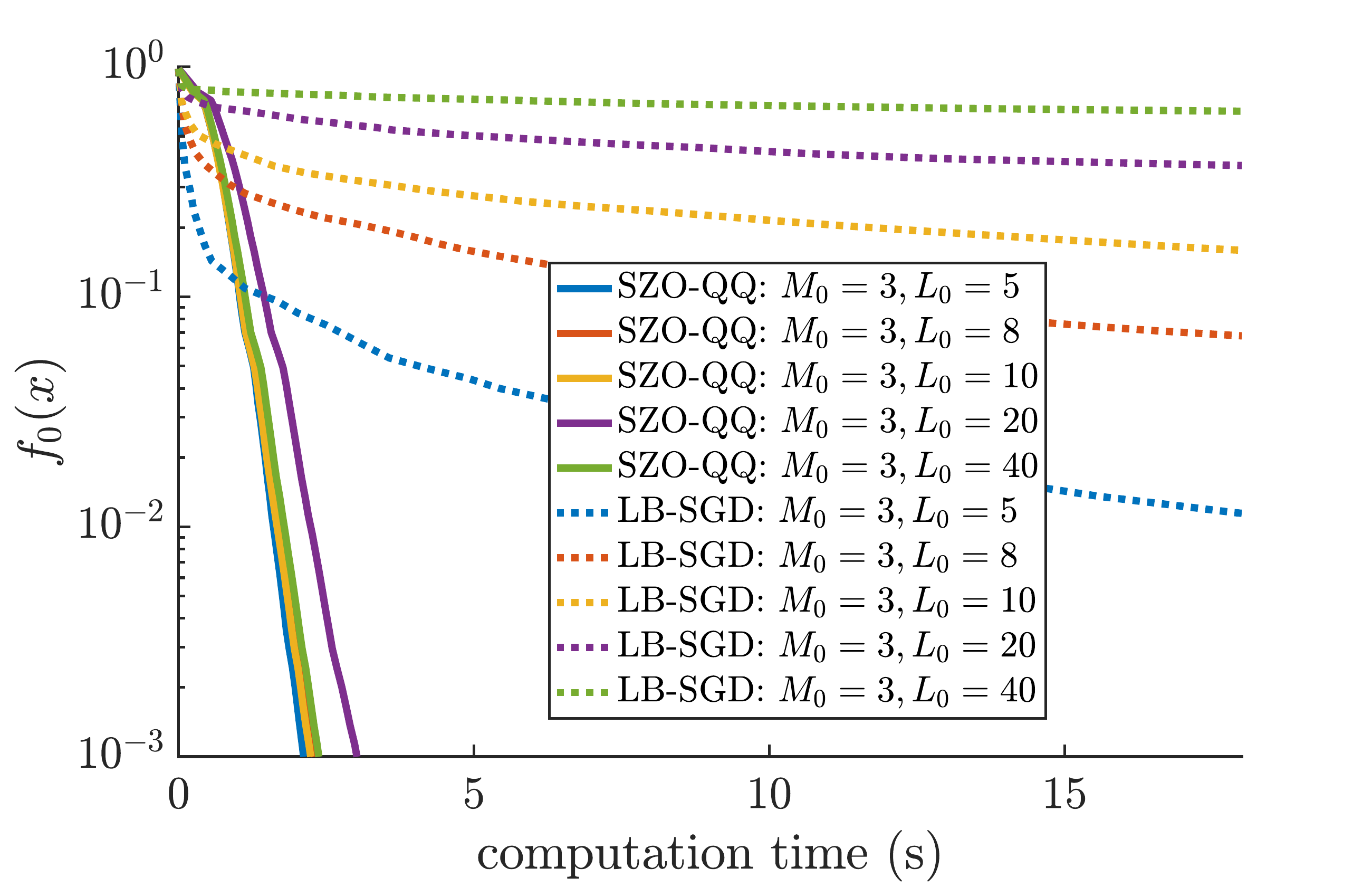}
    \caption{{Objective value as a function of computation time: $\text{ }$ $ M_0=3$ fixed and $L_0$ varied}}
    \label{fig:objective_decrease:Lfixed}
\end{figure}

To show the impact of conservative guesses of $L_i$ and $M_i$, we consider 9 test cases of different values for the pair $(L,M)$. We use $L$ as the Lipschitz constant for all the objective and constraint functions and $M$ as the smoothness constant. We illustrate in Figures \ref{fig:objective_decrease:M fixed} and \ref{fig:objective_decrease:Lfixed} the decrease of the objective function values when SZO-QQ and LB-SGD are applied to solve the 9 test cases. From the figures, we see that the time required by SZO-QQ to achieve an objective function value less than $10^{-2}$ grows with $M$. Despite this, across all the cases SZO-QQ is the first to achieve an objective function value of $10^{-2}$. Another observation is that the performance of SZO-QQ is more sensitive to varying $M$ while LB-SGD is more sensitive to varying $L$. This is due to the differences in the local feasible set formulations in both methods. Indeed, in SZO-QQ the constant $L_i$ is only related to the gradient estimation, and the size of the local feasible set $\mathcal{S}^{(k)}(\cdot)$ in \eqref{eq:safe_set} is mainly decided by $M_i$, while in LB-SGD the size of $\mathcal{T}^{(k)}(\cdot)$, for any $k$, is mainly dependent on the Lipschitz constants $L_i$ for $i\geq 1$. 

We also study the case where the initial guesses of Lipschitz and smoothness constants are wrong, i.e., \eqref{eq: smoothness} in Assumption \ref{ass: smoothness} is violated. With $L=0.2$ and $M=0.2$, we encounter an infeasible sample. Then we follow the method in Remark \ref{rmk: LM} to multiply the constants by 2 every time an infeasible sample is generated. With $L=M=0.8$, every sample is feasible and we derive in 2 seconds an objective function value of $4\times10^{-7}$. In total, we generate two infeasible samples. Although the setting $L=M=0.8$ still fails to satisfy Assumption~\ref{ass: smoothness}, with these constants, SZO-QQ is able to generate iterates that have a subsequence converging to a KKT pair. The readers can check that Theorem \ref{thm: accumulationKKT} holds even if the guesses for Lipschitz and smoothness constants do not verify 
\eqref{eq: smoothness} in Assumption \ref{ass: smoothness} (see the proof of Theorem \ref{thm: accumulationKKT} in Appendix \ref{sec: proof_of_accumulationKKT}).


\subsection{Open-loop optimal control with unmodelled disturbance}
SZO-QQ can be applied to deterministic optimal control problems with unknown nonlinear dynamics by using only feasible samples. To illustrate, we consider a nonlinear system with dynamics $x_{k+1} = Ax_k+Bu_k + \delta(x_k)[1\;0]^\top$, where $x_k\in \mathbb{R}^2$ for $k\geq 0$ and $x_0 = [1\;1]^\top$. The matrices
$$A = \begin{bmatrix}1.1 & 1\\-0.5& 1.1 \end{bmatrix}, B  = \begin{bmatrix}1 & 0 \\0& 1 \end{bmatrix}$$
and the expression of the disturbance $\delta(x):=0.1*(x^{(2)})^2$ are unknown. {The states $\{x_k\}_{k\geq 1}$ solely depend on $\{u_k\}_{k\geq 0}$. We aim to design the inputs $\{u_k\}^5_{k=0}$ (i.e., the decision variables) to minimize the cost
$\sum^{5}_{k=0}\left(x_{k+1}^\top Q x_{k+1} +u_k^\top R u_k\right) $ where $Q= 0.5\,\mathbf I_2$ and $R = 2\,\mathbf I_2$ with identity matrix $\mathbf I_2\in\mathbb R^{2\times 2}$ while enforcing $\|x_{k+1}\|_\infty\leq 0.7$ and $\|u_k\|_\infty\leq 1.5$ for $0\leq k\leq 5$.} Since we assume all the states are measured, we can evaluate the objective and constraints. In this example, we assume to have a feasible sequence of inputs $\{u_k\}_{0\leq k \leq 5}$ that leads to a safe trajectory and results in a cost of 6.81 (as in Assumption \ref{ass: strict feasible}). Different from the settings in the model-based safe learning control methods \cite{fan2020deep,hewing2019cautious}, we do not assume that this safe trajectory is sufficient for identifying the system dynamics with small error bounds. If the error bounds are huge, the robust control problems formulated in \cite{fan2020deep,hewing2019cautious} may become infeasible.

We run SZO-QQ to further decrease the cost resulting from the initial safe trajectory and derive within 146 seconds of computation an input sequence that satisfies all the constraints and achieves a cost of 5.96. This cost is the same as the one obtained when assuming the dynamics are known and applying the solver IPOPT~\cite{wachter2006implementation}. This observation is consistent with Theorem \ref{thm: accumulationKKT} on the convergence to a KKT pair. In this experiment, we set $L_i=M_i=20$ for $i\geq 1$, $\mu = 10^{-4}$ and $\eta = 10^{-1}$. Thus, the parameter $\xi$ adopted is $2\times 10^{-5}$ according to \eqref{eq: xi}.

\subsection{Optimal power flow for an unmodelled electric network}
\label{sec: opf}
{In this section, we apply SZO-QQ to solve the AC Optimal Power Flow (OPF) problem for the IEEE 30-bus system. In this power network, the loads are assumed to be fixed. Six generators are installed, among which one is the slack bus. The slack bus is assumed to provide the active power that is needed to maintain the AC frequency. The 11 decision variables are the voltage magnitude of the 6 generator nodes and the active power generations of the generators, except the slack bus. 142 non-convex constraints are enforced such that the power transmitted through any line is less than the rated value and the voltage magnitude of any bus is within the safe range. The aim of OPF is to minimize a cost, which is a quadratic function of the power generations, such that all the constraints are satisfied. {The ground-truth formulations of the objective and constraint functions can be found in \cite[Section 5.1]{guo2023safe}.}}

{In real-world power systems, it is often hard to formulate the objective and constraint functions of the OPF problem due to unmodelled disturbances (including aging of the devices and external attacks \cite{chuMitigating22}). In \cite{arnold2015model}, extremum seeking is used to solve OPF in a model-free way. As discussed in Section \ref{sec: toy example}, using perturbation signals makes it difficult to select the barrier function coefficient. Considering this issue, we solve OPF by applying SZO-QQ which is also model-free but does not rely on barrier functions. Given a set of values for all the 11 decision variables, we can utilize \texttt{Matpower}~\cite{zimmerman2010matpower} for simulating the power system to 
sample the voltages of all the 30 buses and the power flow through all the transmission lines in this network. Therefore, given any fixed values of the decision variables, we have access to the values of the objective and constraint functions. We also assume to have initial values for all the decision variables such that the constraints are satisfied. In practice, initial values of the decision variables verifying the safety constraints in power systems are not hard to find due to various mechanisms for robust operation, e.g., droop control for power generation, shunt capacitor control and load shedding. }
\begin{figure}[htbp!]
    \centering
    \includegraphics[width=0.7\textwidth]{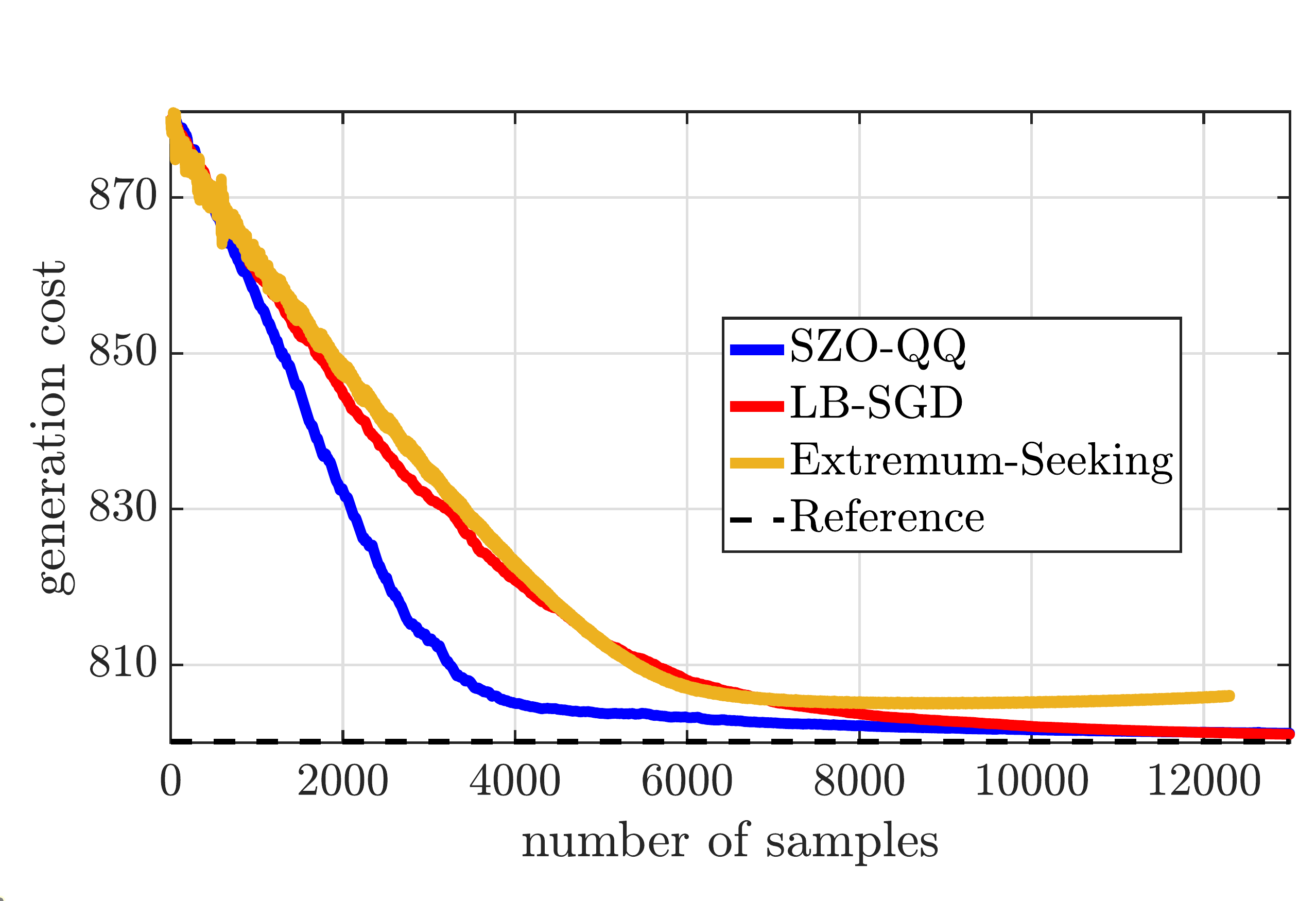}
    \caption{Objective value as a function of the number of samples}
    \label{fig:objective_decrease:opf}
\end{figure}

In this experiment, we set $\mu=0.001$, $\xi = 0.002$, $\Lambda=2$, $M_i=0.2$ and $L_i=1$ for $i\geq 1$. {In Figure \ref{fig:objective_decrease:opf}, we illustrate the decrease of the cost versus the number of samples. For applications to power systems, the time spent on sampling is way more than that on computation. This is because the set points for generators can only vary slowly due to the non-trivial time constants of frequency dynamics. The reference in Fig. \ref{fig:objective_decrease:opf} is derived by assuming to know the ground-truth optimization problem formulation and using Gurobi~\cite{gurobi} to solve it. We see that the achieved cost within 10000 samples is close to the reference. From the comparison with LB-SGD and Extremum-Seeking, we see that to achieve a cost function of 810, SZO-QQ requires around 3200 samples while Extremum-Seeking and LB-SGD need more than 5300, which shows that SZO-QQ is more sample-efficient. Meanwhile, we also compare our method with SafeOptSwarm, which achieves a cost of 823 by using only around 100 samples. However, the cost does not continue to decrease afterwards. We terminate the experiment using SafeOptSwarm after the 1000-th sample, which takes around 9 hours' computation compared to 280 seconds required by SZO-QQ within 4800 samples. From these observations, we see that using SafeOptSwarm to find an initial point for SZO-QQ might be beneficial when sampling takes way more time than computation.}

{We also notice that SZO-QQ and LB-SGD satisfy sample feasibility while Extremum-Seeking generates infeasible iterates. Due to the lack of a recovering mechanism, after the first infeasible sample, most of the other samples of Extremum-Seeking are also infeasible. We discontinue the corresponding curve in Fig. \ref{fig:objective_decrease:opf} after the first infeasible sample. For SZO-QQ, in case an infeasible iterate appears, which is a result of underestimated Lipschitz and smoothness constants, we can easily recover feasibility by using the last feasible iterate. Then, we need to enlarge the constants according to Remark \ref{rmk: LM}.}

\section{Conclusions}
For safe black-box optimization problems, we proposed a method, SZO-QQ, based on samples of the objective and constraint functions to iteratively optimize over local feasible sets. Each iteration of our algorithm involves a QCQP subproblem, which can be solved efficiently. We showed that a subsequence of the algorithm's iterates converges to a KKT pair and no infeasible samples are generated. Given any $\eta>0$, we proposed termination conditions dependent on $\eta$ such that the values returned by the algorithm form an $\eta$-KKT pair. The number of samples required is shown to be $O(\frac{d^2}{\eta^2})$ under mild assumptions. In comparison, the state-of-the-art methods LB-SGD and SafeOpt require $O(\frac{d^2}{\eta^5})$ and $O(\frac{\rho(d)}{\eta^2})$ samples respectively, where $\rho(d)$ is related to the discretization of the space. From numerical experiments, we see that our method can be faster than existing zeroth-order approaches, including LB-SGD, SafeOptSwarm, and Extremum Seeking. Furthermore, the results derived by SZO-QQ are very close to those generated based on the ground-truth model. Future research directions include deriving a tighter complexity bound and generalizing the method for guaranteeing safety even when the samples are noisy.

\bibliographystyle{IEEEtran}
\bibliography{IEEEabrv,reference}
\vspace{.5cm}
\textbf{Appendix}
\appendix
\section{Proof of Theorem \ref{thm: safesetconstruction}}
\label{sec: proofofsafeset}
We notice that $\nu_0^*\leq  l_0^*$ and, from Assumption \ref{ass: smoothness}, $f_i(x_0+\nu_0^*e_j)<f_i(x_0)+l_0^*L_{\max}=0$ for any $i\in \mathbb{Z}_1^m, \, j\in \mathbb{Z}_1^d$. which shows the samples' feasibility. {Furthermore, for any $i\in \mathbb{Z}^m_1$, we have that with $x=x_0$
$$f_i(x_0)+{\nabla}^{\nu_0^*} f_i\left(x_0\right)^\top(x-x_0)+2M_i\|x-x_0\|^2=f_i(x_0)<0,$$
which is to say that $x_0$ belongs to $S^{(0)}(x_0)$.}
To show the feasibility of $\mathcal{S}^{(0)}(x_0)$,
we first partition $\mathcal{S}^{(0)}_i(x_0)$ as 
\[
\mathcal{S}^{(0)}_i(x_0) = \left(\mathcal{S}^{(0)}_i(x_0)\bigcap \mathcal B_{l^*_{0}}(x_0)\right)\left.\left.\bigcup\right(\mathcal{S}^{(0)}_i(x_0)\setminus \mathcal B_{l^*_{0}}(x_0)\right)
\] 
and notice that $\mathcal{S}^{(0)}_i(x_0)\bigcap \mathcal B_{l^*_{0}}(x_0)\subseteq \Omega$. Then, it only remains to show $\mathcal{S}^{(0)}_i(x_0)\setminus \mathcal B_{l^*_{0}}(x_0)\subseteq \Omega.$

For $x\in \mathcal{S}^{(0)}_i(x_0)\setminus B_{l^*_{0}}(x_0)$, we have ${\sqrt{d}\nu^*_0}= l^*_{0} \leq\|x-x_0\|.$
By the mean value theorem, for any $i$ there exists $\theta_i\in [0,1]$ such that 
\begin{align}\label{eq: meanvaluetheorem}
f_i(x)
=&  f_i(x_0)+\nabla f_i\left(x_0 + \theta_i(x-x_0)\right)^\top(x-x_0)\\\notag
= &f_i(x_0) + \nabla f_i\left(x_0 \right)^\top(x-x_0) +\\\notag
& \left(\nabla f_i\left(x_0 + \theta_i(x-x_0)\right) - \nabla f_i\left(x_0 \right)\right)^\top(x-x_0)\\\notag
<& f_i(x_0) + {\nabla}^{\nu^*_0} f_i\left(x_0\right)^\top(x-x_0)\\\notag
&+ \frac{\sqrt{d}\nu^*_0 M_i}{2}\|x-x_0\| + M_i\|x-x_0\|^2\\\notag
\leq& f_i(x_0) + {\nabla}^{\nu^*_0} f_i\left(x_0\right)^\top(x-x_0) + 2M_i\|x-x_0\|^2\\\notag
\leq& 0,
\end{align}
where the first inequality is due to Assumption \ref{ass: smoothness} while the second one can be derived from the definition of $\nu^*_0$ in Theorem \ref{thm: safesetconstruction}. Hence, $\mathcal{S}^{(0)}_i(x_0)\setminus \mathcal B_{l^*_{0}}(x_0)\subseteq \Omega.$ Since $\mathcal{S}^{(0)}_i(x_0)\subseteq \Omega, \forall i$, then $\mathcal{S}^{(0)}(x_0)\subseteq \Omega$. Considering that \eqref{eq: meanvaluetheorem} is a strict inequality, any element of $\mathcal{S}^{(0)}(x_0)$ is strictly feasible for \eqref{eq: optimization problem}. 

\section{Comparison between two different formulations of local safe sets}
\label{sec: comparison of local safe sets}
The works ~\cite{usmanova2022log,sui2015safe} adopt an alternative approximation of the constraints and, in particular, form a local feasible set 
\[
\mathcal{T}^{(0)}(x_0):=\bigcap^m_{i=1}\left\{x:\|x-x_0\|\leq -\frac{f_i(x_0)}{L_i}\right\}.
\]
We see that $\mathcal{T}^{(0)}(x_0)= \{x:F^L_i(x)\leq0,\forall i \}$ where $F^L_i(x):=f_i(x_0)+L_i\|x-x_0\|$ is linear in $\|x-x_0\|$. In contrast, $\mathcal{S}^{(0)}(x_0)=\{x:F^M_i(x)\leq0,\forall i \}$ where $F^M_i(x):=f_i(x_0)+{\nabla}^{\nu^*_0} f_i\left(x_0\right)^\top(x-x_0)+2{M}_i\|x-x_0\|^2$ is a quadratic approximation of $f_i(x)$. In the following proposition, we show that if $x_0$ is sufficiently close to the boundary of the feasible set, $\mathcal{T}^{(0)}(x_0)\subset \mathcal{S}^{(0)}(x_0)$, which means that $\mathcal{S}^{(0)}(x_0)$ is less conservative.
\begin{proposition}
\label{prop: conservative of L}
Let ${\ell_{\min}} = \min_{i\geq 1} (L_i-L_{i,\mathrm{inf}})$. For $x_0$, if 
\begin{equation}
\label{eq: conditionforconservativeness}
{\min_{i\geq 1}} -f_i(x_0)\leq \frac{L_{\max}{\ell_{\min}}}{4{M_{\max}}}, 
\end{equation}
then $\mathcal{T}^{(0)}(x_0)\subset \mathcal{S}^{(0)}(x_0)$.
\end{proposition}
\textbf{Proof.} For any $x\in \mathcal{T}^{(0)}(x_0)$, we have
$$\|x -x_0\|\leq {\frac{\min -f_i(x_0)}{L_{\max}}}\leq{\frac{{\ell_{\min}}}{4M_{\max}} }$$
and thus
\begin{equation}
\label{eq: second_order_inequality}
2{M_{\max}}\|x -x_0\|^2\leq \frac{\ell_{\min}}{2}\|x-x_0\|.
\end{equation}
Considering $\nu^*_0 = \frac{\min-f_i(x_0)}{\sqrt{d}{{L_{\max}}}}$ and \eqref{eq: gradient approximation error}, we have
$$\left\|\Delta^{\nu^*_0}_i(x)\right\|_{2} \leq \frac{\ell_{\min}}{2}$$
and thus 
\begin{equation}
\label{eq: first_order_inequality}
\begin{aligned}
\nabla^{\nu^*_0}f_i(x_0)^\top(x-x_0)=&(\nabla f_i(x_0)+\Delta_i^{\nu^*_0}(x))^\top(x-x_0)\\
\leq&L_{i,\mathrm{inf}} \|x-x_0\|+\frac{\ell_{\min}}{2}\|x-x_0\|
\end{aligned}
\end{equation}
By summing up \eqref{eq: second_order_inequality} and \eqref{eq: first_order_inequality}, we have for any $x\in \mathcal{T}^{(0)}(x_0)$
$$F^M_i(x) \leq f_i(x_0)+ {{L_i}}\|x-x_0\|=F^L_i(x)\leq 0,$$
and thus $x\in \mathcal{S}^{(0)}(x_0)$. \hfill$\blacksquare$
{
\section{Proof of Proposition \ref{prop:feasibility of (SP)}}
\label{prf: feasibility of sp}
Recall From Assumption \ref{ass: strict feasible} that $x_0$ is strictly feasible. Then from Theorem \ref{thm: safesetconstruction} and its proof, we know that $x_1$ is strictly feasible with respect to \eqref{eq: optimization problem}. 
By applying Theorem \ref{thm: safesetconstruction} again through replacing $x_0,l^*_0,\nu^*_0$ with $x_1,l^*_1,\nu^*_1$, we have that for $k = 1$ the iterate $x_1$ is strictly feasible w.r.t. (SP1) and $x_2$ is strictly feasible w.r.t. \eqref{eq: optimization problem}. Through induction, we show the proposition. \hfill $\blacksquare$}

\section{Proof of Proposition \ref{thm: main_theorem}}
\label{sec: 1st argument of thm2}
Proof of Point 1. Given any $k$, we have $x_k\in \mathcal{S}^{(k)}(x_k)$ and $x_{k+1} = \argmin_{x\in \mathcal{S}^{(k)}(x_{k})} f_0(x)+\mu \|x-x_k\|^2$. Thus, 
\begin{equation}
\label{eq:imp_inequality}
\begin{aligned}
&f_0(x_{k+1})+\mu\|x_{k+1}-x_{k}\|^2 \\
\leq\;&  f_0(x_{k})+\mu\|x_{k}-x_{k}\|^2= f_0(x_k).
\end{aligned}
\end{equation}
Proof of Point 2. For $k\geq 0$, one has $f_0(x_k)\leq f_0(x_0)< \beta$ according to Assumption~\ref{ass: bounded_sublevel}. Now we know $\{x_k\}_{k\geq 1}$ is within the set $\mathcal{P}_\beta$. Due to the boundedness of the set $\mathcal{P}_\beta$, we can use Bolzano–Weierstrass theorem to conclude that there exists a subsequence of $\{x_k\}_{k\geq 1}$ that converges. Hence, $\{x_k\}_{k\geq 1}$ has at least one accumulation point $x_c$. According to \eqref{eq:imp_inequality}, $f_0(x_{k+1})\leq  f_0(x_1)-\mu\sum_{i=1}^k \|x_{i+1}-x_{i}\|^2$. Since $f_0(x)$ is a continuous function on the compact set $\mathcal{P}_\beta$, $\mathrm{inf}_{x\in\mathcal{P}_\beta} f_0(x)>-\infty$. Therefore, $\sum_{i=1}^\infty \|x_{i+1}-x_{i}\|^2<\infty$ and $\|x_{k+1}-x_{k}\|$ converges to 0.

Proof of Point 3. The sequence $\{f_0(x_k)\}_{k\geq 1}$ converging to $f_0(x_c)$ is a direct consequence of Point 2 in Proposition \ref{thm: main_theorem} and the continuity of $f_0(x)$. 

\section{Preliminary results towards the proof of Theorem \ref{thm: accumulationKKT}}
\label{sec: preliminary results for Lemma xcandlambdac}
In this section, we only consider the case where $\mathcal{A}(x_c)\neq \emptyset$. Before stating the preliminary results, we have the following notations on the local feasible set $\mathcal{S}^{(k)}(x_k)$ of \eqref{eq: starstar} in Algorithm \ref{alg:sslo}. Our preliminary results are on the properties of the ``limit'' of these feasible sets as $k$ goes to infinity. We define for strictly feasible $x\in \Omega$,
\begin{equation}
\label{eq: center and radius}
\begin{aligned}
O^{(k)}_i(x) :=& x-\frac{\nabla^{\nu^*_{k}} f_i(x)}{2M_i},\\ 
\left(r^{(k)}_i(x)\right)^2 :=&-\frac{f_i(x)}{M_i}+\frac{\|\nabla^{\nu^*_{k}} f_i(x)\|^2}{4M_i^2},
\end{aligned}
\end{equation} 
which allows us to write $$\mathcal{S}^{(k)}_i(x_k) = \mathcal{B}_{r^{(k)}_i(x_k)}(O^{(k)}_i(x_k)).$$
We let $\{x_{k_p}\}_{p\geq 1}$ be a subsequence converging to $x_c$. Since $\{\nu^*_k\}_{k\geq 1}$ converges to 0 (see \eqref{eq: stepsize_2}), we have 
\begin{equation}
\label{eq: convergence_of_origin_radius}
O^{(k_p)}_i(x_{k_p})\rightarrow O_i(x_c),\;r^{(k_p)}_i(x_{k_p})\rightarrow r_i(x_c)\;\text{as}\,p\rightarrow \infty,
\end{equation}
where 
\[
O_i(x_c) := x_c-\frac{\nabla f_i(x_c)}{2M_i},\; \left(r_i(x_c)\right)^2 :=-\frac{f_i(x)}{M_i}+\frac{\|\nabla f_i(x_c)\|^2}{4M_i^2}.
\]
Then, we write
\[
\mathcal{S}_i(x_c) := \mathcal{B}_{r_i(x_c)}(O_i(x_c)),\mathcal{S}(x_c):=\bigcap^m_{i=1}\mathcal{S}_i(x_c).
\]
With these notations, we can state and prove the following lemmas on the properties of $\mathcal{S}(x_c)$.
\begin{lemma}
\label{lmm: preliminary results for Lemma xcandlambdac}
Let Assumption \ref{ass:LICQ} holds and $\mathcal{A}(x_c)\neq \emptyset$ hold, where $\mathcal{A}(x_c): = \{i: f_i(x_c)=0\}$, then
\begin{enumerate}
    \item there exists $x\in\Omega$ that is strictly feasible with respect to $\mathcal{S}(x_c)$, i.e. for any $i\geq 1$, $f^c_i(x)<0$ where
\begin{equation*}
f^c_i(x):=f_i(x_c)+\nabla f_i(x_c)^\top (x-x_c)+2M_i\|x-x_c\|^2.
\end{equation*}
For any $x\in \{x:f^c_i(x)<0, \forall i\}$, there exists $k\in \mathbb{N}$ such that $x$ belongs to $\mathcal{S}^{(k)}(x_k)$;
\item there exists $x_s\in \mathcal{S}(x_c)$ such that $x_s\neq x_c$. For any such $x_s$ and any $0<t<1$, we let $x(t) = tx_c+(1-t)x_s$ and have that $x(t)$ is strictly feasible with respect to $\mathcal{S}(x_c)$. 
\end{enumerate}
\end{lemma}
\textbf{Proof of Point 1}. We let $\mathcal{A}(x_c)=\{i_1,\ldots,i_l\}$. There exists $y\in \mathbb{R}^d$ such that 
\begin{equation}
\label{eq: steering}
Jy = \begin{bmatrix} -1\\ \vdots \\ -1\end{bmatrix}\text{, where } J=\begin{bmatrix}\nabla f_{i_1}(x_c)^\top \\
\vdots \\ \nabla f_{i_l}(x_c)^\top\end{bmatrix},
\end{equation}
because $J$ is full row rank due to LICQ. For any $y$ satisfying~\eqref{eq: steering}, if $\epsilon_0 = 1/(4M_{\max}\|y\|)$, then, for any $\epsilon<\epsilon_0$, $x = x_c+\epsilon y/\|y\|$ and $i\in \mathcal{A}(x_c),$
\begin{equation*}
f^c_i(x)=0-\epsilon/\|y\|+2M_i\epsilon^2< 0.
\end{equation*}
Since $f_i(x_c)<0$ for any $i\notin \mathcal{A}(x_c)$, there exists $\epsilon_c>0$ such that, for $\epsilon<\epsilon_c$ and $x = x_c+\epsilon_c y/\|y\|$, $f^c_i(x)<0$ for any $i\notin \mathcal{A}(x_c)$. Thus, with $\epsilon = \min\{\epsilon_0,\epsilon_c\}$ and $x = x_c+\epsilon y/\|y\|$,  we have $f^c_i(x)<0$ for any $i$. Since $x_c$ is an accumulation point, there exists $k\in \mathbb{N}$ such that $x$ belongs to $\mathcal{S}^{(k)}(x_k)$. 

\textbf{Proof of Point 2}. We utilize the first point and the fact that $x_c$ is not strictly feasible with respect to $\Omega$ to conclude that $x_c$ is not strictly feasible either with respect to $\mathcal{S}(x_c)$ and thus there exists $x_s\in \mathcal{S}(x_c)$ verifying $x_s\neq x_c$. Considering that $f^c_i(x)$ is strongly convex, we have, for any $i$ and any $0<t<1$, $ f^c_i(x(t))<\max\{f^c_i(x_c),f^c_i(x_s)\}\leq 0. $ 

\begin{lemma}
\label{lmm: non_empty_interior}
Let Assumption \ref{ass:LICQ} and $\mathcal{A}(x_c)\neq \emptyset$ hold, then $x_c$ is the unique optimum of the convex optimization 
\begin{equation}
\label{eq: x_c_optimization}
\min_{x\in \mathcal{S}(x_c) } f_0(x)+\mu\|x-x_c\|^2.
\end{equation}
Moreover, the optimizer $\lambda_c$ for the dual variable of \eqref{eq: x_c_optimization} is also unique.
\end{lemma}
\textbf {Proof.} We prove the optimality of $x_c$ by contradiction. Assume $x_c$ is not the optimum of \eqref{eq: x_c_optimization}, then there exists $x_s\in \mathcal{S}(x_c)$ verifying $f_0(x_s)+\mu\|x_s-x_c\|^2< f_0(x_c)$. According to the second point of Lemma \ref{lmm: preliminary results for Lemma xcandlambdac} in Appendix \ref{sec: preliminary results for Lemma xcandlambdac} and the continuity of $f_0(x)$, there exists $0<t<1$ such that with $x(t) = tx_c+(1-t)x_s$ we have $f_0(x(t))+\mu\|x(t)-x_c\|^2< f_0(x_c)$ and $x(t)$ is strictly feasible with respect to $\mathcal{S}(x_c)$. We let $\{x_{k_p}\}_{p\geq 1}$ be a subsequence of $\{x_k\}_{k\geq 1}$ that converges to $x_c$. Considering the first point of Lemma \ref{lmm: preliminary results for Lemma xcandlambdac}, there exists $p$ such that $x(t)\in \mathcal{S}^{(k_p)}(x_{k_p})$. Because of the convergence of the subsequence, we can assume without loss of generality that $p$ is sufficiently large so that $f_0(x(t))+\mu\|x(t)-x_{k_p}\|^2< f_0(x_c)$. Due to the optimality of $x_{k_p+1}$ for the problem \eqref{eq: starstar} in Algorithm~\ref{alg:sslo} when $k=k_p$, we have $f_0(x_{k_p+1})+ \mu\|x_{k_p+1}-x_{k_p}\|^2<f_0(x_c)$, which contradicts the monotonicity of $\{f_0(x_k)\}_{k\geq 1}$ in Proposition \ref{thm: main_theorem}.  Due to optimality of $x_c$ and LICQ, there exists $\lambda_c\in \mathbb{R}^m$ such that $(x_c,\lambda_c)$ is a KKT pair of \eqref{eq: x_c_optimization}.

We prove the uniqueness of $x_c$ and $\lambda_c$ also by contradiction. Assume there exists $x_o\in \mathcal{S}(x_c)\setminus\{x_c\}$ such that  $f_0(x_o)+\mu\|x_o-x_c\|^2= f_0(x_c)$. Due to the strong convexity of the function $f_0(x)+\mu\|x-x_c\|^2$, we know
$$f_0(\frac{x_o+x_c}{2})+\mu\|\frac{x_o+x_c}{2}-x_c\|^2< f_0(x_c),$$
which contradicts the optimality of $x_c$ for \eqref{eq: x_c_optimization}. Assume $(x_c,\lambda_{c,1})$ and $(x_c,\lambda_{c,2})$ are two KKT pairs of \eqref{eq: x_c_optimization} with $\lambda_{c,1}\neq \lambda_{c,2}$, then for $j=1,2$, $\lambda_{c,j}^{(i)}=0$ for any $i\neq \mathcal{A}(x_c)$,
$$ \sum^m_{i=1}\lambda_{c,j}^{(i)}\nabla f_i(x_c) = - \nabla f_0(x_c) $$
and thus 
$$ \sum^m_{i=1}(\lambda_{c,1}^{(i)}-\lambda_{c,2}^{(i)})\nabla f_i(x_c)=0.$$
This contradicts LICQ at $x_c$ since $\lambda_{c,1}-\lambda_{c,2}\neq 0.$ \hfill$\blacksquare$

\section{Proof of Theorem \ref{thm: accumulationKKT}}
\label{sec: proof_of_accumulationKKT}
This proof only considers the case where $\mathcal{A}(x_c)\neq \emptyset$ and can be easily extended to ``$\mathcal{A}(x_c)= \emptyset$''. According to Lemma \ref{lmm: non_empty_interior}, there exists a $\lambda_c\in\mathbb{R}^m$ such that $(x_c,\lambda_c)$ is a KKT pair of \eqref{eq: x_c_optimization}, i.e., 
$$
\begin{aligned}
&\nabla f_0(x_c) + \sum_{i\in\mathcal{A}(x_c)}\lambda_c^{(i)}\nabla f_i(x_c) = 0 \\
\text{and }\; &\lambda_c^{(i)} = 0\text{ for }i\notin\mathcal{A}(x_c)\text{, } 
\end{aligned}
$$ 
which coincides with KKT conditions of \eqref{eq: optimization problem}. Thus, $(x_c,\lambda_c)$ is also a KKT pair of of \eqref{eq: optimization problem}. Following the same arguments used in the proof of Lemma \ref{lmm: non_empty_interior}, one can exploit LICQ to show that there does not exist $\lambda_{c,2}\neq \lambda_c$ such that $(x_c,\lambda_{c,2})$ is a KKT pair of \eqref{eq: optimization problem}. 

\section{Geometric illustration of an upperbound to $\|\lambda_c\|_\infty$}
\label{sec: characterization of lambda}
\begin{figure}[htbp!]
    \centering
\includegraphics[width=0.5\linewidth]{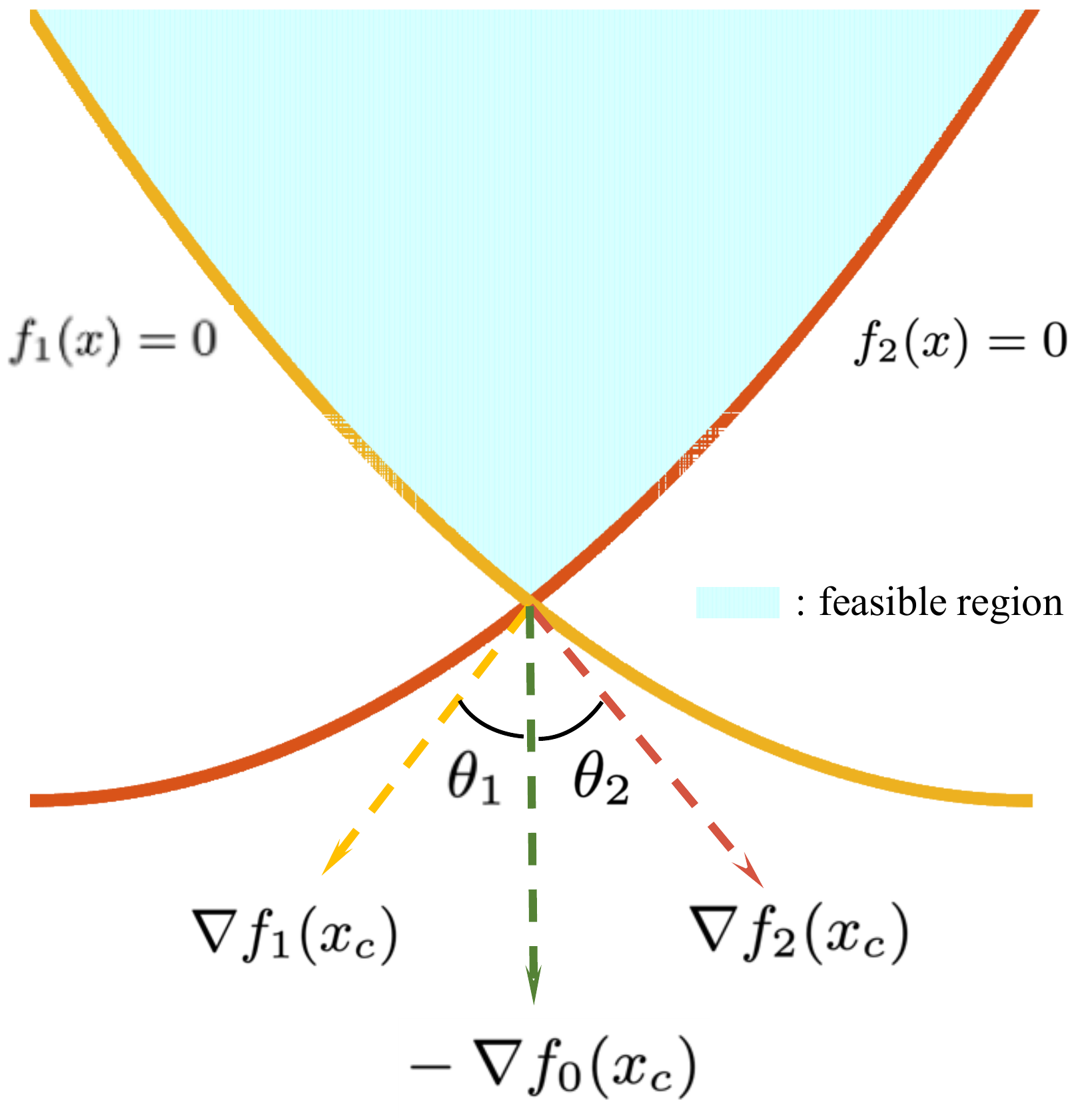}
    \caption{The illustration of angles $\theta_1$ and $\theta_2$}
    \label{fig:illustration}
\end{figure}

With the following example, we aim to illustrate that $\|\lambda_c\|_\infty$ is related to the geometric properties of the feasible region. We consider an instance of the optimization problem \eqref{eq: optimization problem} where $d=2$, the feasible region is convex, and $(x_c,\lambda_c)$ is a KKT pair. We only consider the non-degenerate case where $\mathcal{A}(x_c)=\{1,2\}$ and assume LICQ holds at $x_c$, i.e., $\nabla f_i(x_c)$ are linearly independent for $i=1,2.$ The objective and constraint functions are normalized at $x_c$, i.e., $\|\nabla f_i(x_c)\| = 1$ for $i\in\mathbb{Z}^{2}_{0}.$ Then, we use coordinate transformation such that $\nabla f_0(x_c) = \begin{bmatrix}0 &  -1\end{bmatrix}^\top.$ Since $\mathcal{A}(x_c)=\{1,2\}$, the KKT pair $(x_c,\lambda_c)$ satisfies that $\lambda_c^{(1)},\lambda_c^{(2)}\geq 0$ and
\begin{equation}
\label{eq: example_kkt}
\begin{aligned}
 f_1(x_c)\leq 0, \text{ } f_2(x_c)&\leq 0\\
\nabla f_0(x_c)+\lambda^{(1)}_c\nabla f_1(x_c)+\lambda^{(2)}_c\nabla f_2(x_c) & = 0\\
 \lambda^{(1)}_c f_1(x_c) = 0,\text{ } \lambda^{(2)}_c f_2(x_c) & = 0.
\end{aligned}
\end{equation}
Let $\theta_i$ be the angle between $-\nabla f_0(x_c)$ and $\nabla f_i(x_c)$ for $i = 1,2.$ Due to the convexity of the feasible region, $0<\theta_1+\theta_2<\pi.$ By solving \eqref{eq: example_kkt}, we have that 
$$\lambda^{(1)}_c = \frac{|\sin \theta_2|}{\sin(\theta_1+\theta_2)},\;\lambda^{(2)}_c = \frac{|\sin \theta_1|}{\sin(\theta_1+\theta_2)}.$$
We illustrate in Fig. \ref{fig:illustration} how to construct $\theta_1$ and $\theta_2$.

We notice that  
$$\|\lambda_c\|_\infty<(\sin\theta)^{-1},$$
where $\theta = \pi -\theta_1-\theta_2$ is the angle between the two lines $l_1:=\{x:(x-x_c)^\top\nabla f_1(x_c)=0\}$ and $l_2:=\{x:(x-x_c)^\top\nabla f_2(x_c)=0\}$. These two lines are actually the boundaries formed by the constraint functions $f_1(x)$ and $f_2(x)$ linearized at $x=x_c$. From the above conclusions, we see that for 2-dimensional optimization problems, we need a large $\Lambda$ to satisfy Assumption \ref{ass: magnitude_dual_variable} only when the angle $\theta$ is small.

{It is possible to extend the above analysis to problems with a higher dimension. Then, with $\mathcal{A}(x_c)=\{i_1,\ldots,i_k\}$ and  $F_{\mathcal{A}} := [\nabla f_{i_1}(x_c)\; \nabla f_{i_2}(x_c)\ldots \nabla f_{i_k}(x_c)]^\top$, one can show that $\sigma^{-1}_{\min}(F_\mathcal{A})$ is an  upperbound for $\|\lambda_c\|_\infty$, where $\sigma_{\min}(F_\mathcal{A})$ denotes the smallest singular value of $F_\mathcal{A}$.}

\section{Proof of Proposition \ref{thm: algorithm_termination}}
\label{sec: algorithm_termination}

Since $(x_c,\lambda_c)$ is a KKT pair where $\|\lambda_c\|_{\infty}<\Lambda$, by using triangular inequalities on the norm terms defining $\delta_1(k,\lambda_c)$, we obtain 
\begin{align*}
&\delta_1(k,\lambda_c)\\
\leq & \|{\nabla} f_0(x_c)  + \sum^{m}_{i=1} \lambda_{c}^{(i)}{\nabla} f_i\left(x_{c}\right) \|+ \|{\nabla} f_0(x_k) - {\nabla} f_0(x_c) \| \\
&+4\mu \|x_{k+1}-x_k\|
+\sum^{m}_{i=1} \Lambda\Big(\|{\nabla}^{\nu_{k}^*} f_i\left(x_{k}\right)-{\nabla} f_i\left(x_{k}\right)\| \\&+ \|{\nabla} f_i\left(x_{k}\right)- {\nabla} f_i\left(x_{c}\right)\|\Big)+4m{M_{\max}}\Lambda\|x_{k+1}-x_{k}\|.
\end{align*}
Similar computations for $\delta^{(i)}_2(k,\lambda_c)$ give
\begin{align*}
\delta^{(i)}_2(k,\lambda_c)
\leq&  |\lambda^{(i)}_c f_i(x_c)|+\Lambda|f_i(x_k)-f_i(x_c)| 
\\&+\Lambda\Big( L_i|x_{k+1}-x_k|+2M_i\|x_{k+1}-x_k\|^2\Big).
\end{align*}
We let $\{x_{k_p}\}_{p\geq 1}$ be an subsequence that converges to $x_c$. Considering that the gradient estimation error converges to 0 (see \eqref{eq: stepsize_2} and Lemma \ref{lmm: gradient estimation error}), we know the term $\|{\nabla}^{\nu_{k_p}^*} f_i\left(x_{k_p}\right)-{\nabla} f_i\left(x_{k_p}\right)\|$ converges to 0 as $p$ goes to infinity. Therefore, we have $$\lim_{p\rightarrow \infty}\delta_1(k_p,\lambda_c) = 0\text{  and }\lim_{p\rightarrow \infty}\max_{1\leq i \leq m}|\delta^{(i)}_2(k_p,\lambda_c)| = 0.$$
Thus, for any $k_0$ and $\eta>0$, one can find $k_{\Lambda}>k_0$ such that $\max\{\delta_1(k_{\Lambda},\lambda_c),\max_{1 \leq i\leq m}|\delta^{(i)}_2(k_{\Lambda},\lambda_c)| \}<\eta/2$. For $k=k_{\Lambda}$ in SZO-QQ, we have that $\lambda_{k_{\Lambda}+1}$, the solution to \eqref{eq: SP2}, has an infinite norm less than $2\Lambda$ because $\lambda_c$ is a feasible solution to \eqref{eq: SP2} and $\|\lambda_c\|_\infty<2\Lambda$, which is to say that the second termination condition of Algorithm \ref{alg:sslo} is satisfied when $k=k_{\Lambda}$.

Since $\|x_{k+1}-x_k\|$ converges to 0 as $k$ goes to infinity (see Proposition \ref{thm: main_theorem}), we can choose $k_0$ to be sufficiently large so that, for any $k>k_0$, $\|x_{k+1}-x_k\|\leq \xi$. Thus, the two termination conditions are satisfied when $k=k_{\Lambda}$. 

\section{Proof of Theorem \ref{prop: eta-KKT} }
\label{sec: proof_of_xi}
The pair $(\tilde{x},\tilde{\lambda})$ and the index $\tilde{k}$ returned by Algorithm \ref{alg:sslo} satisfy $\tilde{x} = x_{\tilde{k}}$ and 
\begin{equation}
\label{eq: inequality_delta}
\max\Big\{\delta_1(\tilde{k},\tilde{\lambda}), \max\{\delta^{(i)}_2(\tilde{k},\tilde{\lambda}):i\geq 1 \}\Big\}\leq\frac{\eta}{2}.
\end{equation}
By using triangular inequalities we have for any $i$,
\begin{align}\notag
& \|{\nabla}^{\nu_{\tilde{k}-1}^*} f_i\left(x_{\tilde{k}-1}\right) - {\nabla} f_i\left(x_{\tilde{k}}\right)\|\\\notag
\leq & \|{\nabla}^{\nu_{\tilde{k}-1}^*} f_i\left(x_{\tilde{k}-1}\right) - {\nabla} f_i\left(x_{\tilde{k}-1}\right)\| \\\notag
&\qquad\qquad\qquad\qquad + \|{\nabla} f_i\left(x_{\tilde{k}}\right) - {\nabla} f_i\left(x_{\tilde{k}-1}\right)\|\\\label{eq:first_inequality}
\leq & \alpha \xi + L_i\xi.
\end{align}
Then, based on \eqref{eq: xi}, \eqref{eq: inequality_delta} and \eqref{eq:first_inequality}, we have
\begingroup\makeatletter\def\f@size{8}\check@mathfonts
\begin{align*}
&||\nabla f_0\left(x_{\tilde{k}}\right) +\sum^{m}_{i=1} \tilde{\lambda}^{(i)}{\nabla} f_i\left(x_{\tilde{k}}\right)\| \\
\leq & \Big\|\nabla f_0\left(x_{\tilde{k}}\right) +\sum^{m}_{i=1} \tilde{\lambda}^{(i)}\left({\nabla}^{\nu_{\tilde{k}-1}^*} f_i\left(x_{\tilde{k}-1}\right)+4M_i(x_{\tilde{k}}-x_{\tilde{k}-1})\right)\\
&+2\mu(x_{\tilde{k}}-x_{\tilde{k}-1}) \Big\| + \|2\mu(x_{\tilde{k}}-x_{\tilde{k}-1})\|\\
& +2\Lambda \sum^m_{i=1}\Big(4\left\|M_i(x_{\tilde{k}}-x_{\tilde{k}-1})\right\|\Big) \\&+ 2\Lambda \sum^m_{i=1}\Big(\left\|  {\nabla}^{\nu_{\tilde{k}-1}^*} f_i\left(x_{\tilde{k}-1}\right) - {\nabla} f_i\left(x_{\tilde{k}}\right)\right\|\Big)\\
\leq & \eta/2 + 2\Lambda\sum^m_{i=1}\left(5M_i\xi+\alpha_i \nu^*_{\tilde{k}-1}\right) +2\mu\xi \leq  \eta,
\end{align*}
\endgroup
and
\begin{align*}
&\left\|\tilde{\lambda}^{(i)}f_i(x_{\tilde{k}})\right\| \\
\leq & \Big\|\tilde{\lambda}^{(i)} \Big(f_i(x_{\tilde{k}-1})+{\nabla}^{\nu_{\tilde{k}-1}^*} f_i\left(x_{\tilde{k}-1}\right)(x_{\tilde{k}}-x_{\tilde{k}-1})\\
& +2M_i\|x_{\tilde{k}}-x_{\tilde{k}-1}\|^2\Big)\Big\| + 
2\Lambda\Big(\|f_i(x_{\tilde{k}}) - f_i(x_{\tilde{k}-1})\|+\\
&\|{\nabla}^{\nu_{\tilde{k}-1}^*} f_i\left(x_{\tilde{k}-1}\right) - {\nabla} f_i\left(x_{\tilde{k}-1}\right)\|\cdot \|x_{\tilde{k}}-x_{\tilde{k}-1}\|+ \\
& \|\nabla f_i(x_{\tilde{k}-1})\|\cdot\|x_{\tilde{k}}-x_{\tilde{k}-1}\| + 2M_i\|x_{\tilde{k}}-x_{\tilde{k}-1}\|^2\Big)\\
\leq & \eta/2 + 2\Lambda( 2L_i \xi + \alpha_i \xi^2   + 2M_i\xi^2)\\
\leq& \eta/2 + 2\Lambda( 2L_i + \alpha_i + 2M_i)\xi \leq \eta,
\end{align*}
which concludes the proof.

\section{Proof of Lemma \ref{lmm: convergence}}
\label{sec: proof of convergence lemma}
We assume $x_\tau$ is an accumulation point of $\{x_k\}_{k\geq 1}$ and is also a strict local minimizer. We show the convergence of $\{x_k\}_{k\geq 1}$ by contradiction. We assume that $\mathcal{C}\setminus \{x_\tau\}\neq \emptyset$, where $\mathcal{C}$ is the set of accumulation points of $\{x_k\}_{k\geq 1}$. Then, there exists $\epsilon>0$ such that $\mathcal{C}\cap (\Omega\setminus\mathcal{B}_\epsilon(x_\tau))\neq \emptyset$ and any $x\in \mathcal{B}_\epsilon(x_\tau)\setminus\{x_\tau\}$ verifies $f_0(x_\tau)<f_0(x)$. Since the sphere $\mathcal{SP}_\epsilon(x_\tau)=\{x:\|x-x_\tau\| = \epsilon\}$ is compact, we let $\sigma = \inf_{x\in\mathcal{SP}_\epsilon(x_\tau)} f_0(x)$ and have $\sigma>f_0(x_\tau)$. Therefore, there exists $k_\alpha>0$ such that $f_0(x_{k_\alpha})<(\sigma+f_0(x_\tau))/2$. 

Since there exists an accumulation point outside $\mathcal{B}_\epsilon(x_\tau)$ and $\{x_{k+1}-x_k\}_{k\geq 1}$ converges to 0, we can find $k_\beta>k_\alpha$ such that $x_{k_\beta}\in \mathcal{B}_\epsilon(x_\tau)$, $x_{k_\beta+1}\notin \mathcal{B}_\epsilon(x_\tau)$ and $\|x_{k_\beta+1}-x_{k_\beta}\|\leq (\sigma-f_0(x_\tau))/(4L_{\max})$. Let $\tilde{x} = \{x: \text{there exists }t\in[0,1] \text{ such that ${x} = tx_{k_\beta}+(1-t)x_{k_\beta+1}$}\}\cap \mathcal{SP}_\epsilon(x_\tau)$, i.e., $\tilde{x}$ is the intersection of $\mathcal{SP}_\epsilon(x_\tau)$ and the line segment between $x_{k_\beta}$ and $x_{k_\beta+1}$. Then, $\|x_{k_\beta}-\tilde{x}\|\leq (\sigma-f_0(x_\tau))/(4L_{\max})$ and thus
\begin{equation*}
\begin{aligned}
f_0(x_{k_\beta})&\geq f_0(\tilde{x})-\frac{\sigma-f_0(x_\tau)}{4}\\
&\geq \sigma -\frac{\sigma-f_0(x_\tau)}{4}\\
&> (\sigma+f_0(x_\tau))/2> f_0(x_{k_\alpha}),
\end{aligned}
\end{equation*}
which contradicts with the monotonicity of $\{f_0(x_k)\}_{k\geq 1 }$ shown in Proposition \ref{thm: main_theorem}. \hfill$\blacksquare$
  
\section{Proof of Lemma \ref{lmm: convergence of the dual variable}}
\label{sec: lambda_conv}
To begin with, we let $\mathcal{D}_\lambda(y,\nu)\subset \mathbb{R}^m$ be the optimal solution set of the dual of the following convex problem:
\begin{equation}
\begin{aligned}
\mathsf{P}(y,\nu): \min_{x\in\mathbb R^d} & \quad f_0(x)+\mu \|x-y\|^2\\
\text{subject to} &  \quad f_i(y) + \Big(\underbrace{\Delta^\nu_i(y) + \nabla f_i(y)}_{\nabla^\nu_i f_i(y) }\Big)^\top(x-y)\\
&\quad+2M_i\|x-y\|^2\leq 0.
\end{aligned}
\end{equation} 
Notice that the problem  $\mathsf{P}(x_k,\nu^*_k)$ coincides with \eqref{eq: starstar} in Algorithm \ref{alg:sslo}. Then, we have
\[
\begin{aligned}
\mathcal{D}_\lambda&(y,\nu) :=\\
& \argmax_{\lambda\geq 0}\min_{x} f_0(x)+\mu \|x-y\|^2
+\sum_{i=1}^m\lambda^{(i)}\Big(f_i(y)\\
&\quad + \Big(\Delta^\nu_i(y) + \nabla f_i(y)\Big)^\top(x-y) +2M_i\|x-y\|^2\Big).
\end{aligned}
\] 
By solving the inner minimization problem which is an unconstrained convex quadratic programming, we know that there exist $p\in\mathbb{R}^{m\times 1}_{\geq 0}$ and $a\in\mathbb{R}_{>0}$, independent of $y$ and $\nu$, such that
\[
\mathcal{D}_\lambda(y,\nu)= \argmax_{\lambda\geq 0}G(\lambda,y,\nu),
\]
where for $\lambda\in\mathbb{R}^m_{\geq0}$, 
\[
G(\lambda,y,\nu):=\frac{\lambda^\top Q(y,\nu)\lambda + q^\top(y,\nu)\lambda+b(y,\nu)}{p^\top \lambda + a},
\] 
the functions $Q(y,\nu)\in\mathbb{R}^{m\times m}, q(y,\nu)\in\mathbb{R}^{m\times 1}, b(y,\nu)\in\mathbb{R}$ are continuous in $(y,\nu)$ and $Q(y,\nu)$ is negative definite. 

From the continuity of $Q(y,\nu), q(y,\nu)$ and $b(y,\nu)$, the function $G(\lambda,y,\nu)$ is continuous in all arguments for $\lambda\geq 0$. Due to the continuity and  the uniqueness of the optimal dual solution $\mathcal{D}_\lambda(x_c,0) = \{\lambda_c\}$ (see Lemma \ref{lmm: non_empty_interior}), we can use perturbation theory~\cite[Proposition 4.4]{bonnans2013perturbation} to conclude that $\mathcal{D}_{\lambda}(y,\nu)$ is upper semicontinuous at $(y,\nu)=(x_c,0)$. 
\begin{definition}
Let $W$ and $V$ be two vector spaces. A multifunction $F: W\rightarrow \mathcal{P}(V)$, where $\mathcal{P}:=\{P: P\subset V\}$, is said to be upper semicontinuous at $w_0$ if for any neighborhood $\mathcal{N}_{V}$ of $F(w_0)$, there exists a neighborhood $\mathcal{N}_W$ of $w_0$ such that the inclusion $F(w)\subset \mathcal{N}_{V}$ holds for any $w\in\mathcal{N}_W$.
\end{definition}
Considering the convergence of $(x_k,\nu^*_k)$ to $(x_c,0)$ and the upper semicontinuity of $\mathcal{D}_{\lambda}(y,\nu)$ at $(y,\nu)=(x_c,0)$, for any $\delta>0$, there exists $k_\delta>0$ such that $\mathcal{D}_\lambda(x_k,\nu^*_k)\subset\mathcal{B}_\delta(\lambda_c)$ for any $k>k_\delta$. Since $\lambda^\circ_{k+1}\in\mathcal{D}_\lambda(x_k,\nu^*_k)\subset\mathcal{B}_\delta(\lambda_c) $, we have $\lambda^\circ_{k+1}\in\mathcal{B}_\delta(\lambda_c)$ for any $k>k_\delta$.

\section{Proof of Theorem \ref{thm: complexity}}
\label{sec: proof_of_complexity}
According to Lemma \ref{lmm: convergence of the dual variable}, there exists $\bar{k}>0$ such that $\|\lambda^\circ_{k+1}\|_\infty\leq 2\Lambda$ for any $k\geq \bar{k}$. Since $\lambda^\circ_{k+1}$ is a feasible solution of \eqref{eq: SP2} in Algorithm \ref{alg:sslo}, $\lambda_{k+1}$, the optimal solution of \eqref{eq: SP2}, also satisfies  $\|\lambda_{k+1}\|_\infty\leq 2\Lambda$. 

Recall the definition of $h(\eta)$ in \eqref{eq: xi}. We let $\bar{\xi} := \inf_{k\leq \bar{k}} \|x_{k+1}-x_k\|$, $\bar{\eta} := \inf\{\eta: h(\eta)\geq \bar{\xi}/2\}$ and consider the case where $\eta<\bar{\eta}$. We first notice that if $\|x_{k+1}-x_k\|\leq h(\eta)$, we have $\|x_{k+1}-x_{k}\|<\bar{\xi}$ and thus $k> \bar{k}$. We then let $\mathcal{K}(\eta) := \max\{k:\|x_{k+1}-x_k\|>h(\eta)\}+1$. Since $\|x_{\mathcal{K}(\eta)+1}-x_{\mathcal{K}(\eta)}\|\leq h(\eta)$, we have that $\mathcal{K}(\eta)>\bar{k}$ and thus $\|\lambda_{k+1}\|_\infty\leq 2\Lambda$, which is equivalent to say that with $k = \mathcal{K}(\eta)$, the two termination conditions in {Algorithm}~\ref{alg:sslo} are satisfied. Then $\tilde{k}$, the iteration number returned by Algorithm~\ref{alg:sslo}, verifies that $\tilde{k}\leq \mathcal{K}(\eta)+1$. According to \eqref{eq: monotonicity_inequality},
\begin{equation}
\begin{aligned}
f_0(x_0)- \inf\{f_0(x): x\in \Omega\}&\geq f_0(x_0)-  f_0(x_{\mathcal{K}(\eta)})\\
&\geq \mu\mathcal{K}(\eta)(h(\eta))^2
\end{aligned}
\end{equation}
and thus $\tilde{k}\leq \mathcal{K}(\eta)+1\leq \overline{\mathcal{K}}(\eta)+1$. {Therefore, since $f_0(x_0)- \inf\{f_0(x): x\in \Omega\}<C_1$, according to the definition of $h(\eta)$ in \eqref{eq: xi} there exists $C_2>0$ such that $\mathcal{K}(\eta)+1\leq C_2(d/\eta^{2})$. Now, we can conclude that the number of iterations is at most $O(d/\eta^2)$. Since in each iteration, $d+1$ samples are taken, the number of samples required by Algrorithm \ref{alg:sslo} is at most $O(d^2/\eta^2)$.}

\end{document}